\documentclass[12 pt]{article}
\usepackage{graphicx}
\usepackage{subfigure}
\usepackage{amsthm}
\newcommand{\hs}{\hskip 0.2cm}
\newcommand{\hp}{\hskip 1.5mm}

\newcommand{\ha}{\theta}
\newcommand{\va}{\Psi}

\def\be{\begin{equation}}
\def\ee{\end{equation}}

\title{\large\bf Lagrange stability for impulsive Duffing equations}
\author{ Jianhua Shen$^1$, Lu Chen$^1$, Xiaoping Yuan$^{2,}$${\thanks{The corresponding author.
{\it E-mail addresses:} chenlu@zju.edu.cn (L. Chen), jianhuashen2013@163.com (J. Shen), xpyuan@fudan.edu.cn (X. Yuan).}}$\\
\small {\it $^1$Department of Mathematics, Hangzhou Normal University}\\
\small {\it Hangzhou,  Zhejiang 310036, PR China}\\
\small {\it $^2$School of Mathematical Sciences, Fudan University}\\
\small {\it Shanghai 200433, PR China}}
\date{}
\begin{document}
 \maketitle

\noindent {\small\bf Abstract} \vskip 0.2cm

{\small  This work discusses the boundedness of solutions for impulsive Duffing equation with time-dependent polynomial potentials. By KAM theorem, we prove that all solutions of the Duffing equation with low regularity in time undergoing suitable impulses are bounded for all time and that there are many (positive Lebesgue measure) quasi-periodic solutions clustering at infinity. This result extends some well-known results on Duffing equations to impulsive Duffing equations. \\

 }
 \vskip 0.4cm
 \noindent{\bf 1. Introduction}\\

Let us begin with the harmonic oscillator (linear spring)
\[\ddot{x}+k^2\, x=0.\]
All solutions of this equation are bounded for $t\in \mathbf{R}$. That is, this equation is Lagrange stable. However, the stability is violated when the linear spring is stressed by an external force, periodic in time. More exactly, there is a unbounded solution to the equation
\[\ddot{x}+k^2\, x=p(t)\]
where the frequency of $p$ is equal to the frequency $k$ of the spring itself.  Now let us consider a nonlinear equation
\[\ddot{x}+x^3=0.\]
 This equation is a Hamiltonian system whose all solutions are periodic in time, thus, it is Lagrange stable, too. An interesting problem is that whether or not
 \[\ddot{x}+x^3=p(t)\] is stable when $p(t)$ is periodic in time.
   Inspired  by the question of Littlewood [10], Morris [13] in 1976  gave a positive answer. Actually, there is a long history for this subject.
 In 1962, Moser [14] proposed to study the boundedness of all solutions for Duffing equation
$$ \ddot{x}+\alpha x^3+\beta x=p(t), $$
where $\alpha>0, \beta\in \mathbf{R}$ are constants, $p(t)$ is a continuous function with period 1.
  Subsequently, Morris's boundedness result was, by Dieckerhoff-Zehnder [3] in 1987, extended to a family of Duffing equations \\
$$\ddot{x}+x^{2n+1}+\sum_{i=0}^{2n} x^{i}p_{i}(t)=0,\hs n\geq 1,\eqno (*)$$
where $p_i(t)(i=0,1,\cdots,2n)$ are 1-periodic and sufficiently smooth functions. The smoothness has been recently relaxed to
  $C^\gamma$-H\"older continuity with $\gamma>1-\frac1n$ in [29]. See[6, 9, 11, 12, 25, 26-28] for more details.

The equation ($*$) can be regarded as  $\ddot{x}+x^{2n+1}=0$  with a  perturbation $\sum_{i=0}^{2n} x^{i}p_{i}(t)$.
The results mentioned as the above show that the nonlinear equation  $\ddot{x}+x^{2n+1}=0$ is Lagrange stable under a periodic perturbation in the equation itself.

 {\it What happens when the nonlinear equation  $\ddot{x}+x^{2n+1}=0$  is subject to both periodic perturbation and an impulse at the same time?}

In the present paper,  we will discuss the boundedness of solutions and the existence of quasi-periodic solutions for the impulsive Duffing equation\\
  $$\left\{\begin{array}{ll}
  \ddot{x}+x^{2n+1}+\sum_{i=0}^{2n} x^{i}p_{i}(t)=0, \hs t\neq t_j,\hs n\geq1,\\
  \bigtriangleup x(t_{j}):=x(t_{j}^+)-x(t_{j}^-)=I_j(x(t_j^-),\dot{x}(t_j^-)),\\
  \bigtriangleup \dot{x}(t_j):=\dot{x}(t_j^+)-\dot{x}(t_j^-)=J_j(x(t_j^-),\dot{x}(t_j^-)),\hs j=\pm1,\pm2,\cdots,\\
   \end{array}\right.  \eqno (1.1)$$\\
  where for $j=\pm1,\pm2,\cdots, \{t_j\} $ is a strictly increasing sequence of real numbers, which is called an impulsive time sequence, $I_j,J_j: \mathbf{R^2}\to \mathbf{R}$ are the sequences of impulsive functions, and where $p_i(t) (0\leq i\leq 2n)$ are 1-periodic functions.

  There are few results on Lagrange stability and the existence of quasi-periodic solutions of impulsive differential equations. However, there are many studies on the existence of periodic solutions of impulsive differential equations. See [4, 5, 16,17, 20,24], for example. We refer the reader to classical monographs [1, 7] for the general theory of impulsive differential equations.

     \noindent{\bf 2. Statement of results}\\

     To formulate our main result we have to introduce some notations and hypotheses.

    Let $\mathbf{R}, \mathbf{C}, \mathbf{N}$ and $\mathbf{Z}$ be the sets of all real numbers, complex numbers, natural numbers and integers, respectively. Denote by $\mathbf{\mathcal{T}}$ the impulsive time sequence $\{t_j\}, j=\pm1,\pm2,\cdots$, and denote by $\mathbf{\mathcal{A}}$ the set of indexes $j$.

    We say that a function $x(t)$ is a solution of Eq. (1.1) if it satisfies (1.1). In Section 3, we will describe in detail the properties of the solution operators for equation (1.1), and explain why all the solutions of (1.1) exist for all $t\in \mathbf{R}$ under appropriate hypotheses on the impulsive functions.

        In the present paper, we assume that the following condition $({\rm \mathbf{H}})$ holds true. \\

 \noindent $({\rm \mathbf{H}})$  There exists a positive integer $k$ such that $\hp 0< t_1<t_2<\cdots<t_{k}<1$, and that $t_j$'s, $I_j(x,y)$'s, $J_j(x,y)$'s are $k$-periodic in $j$ in the sense that $t_{j+k}=t_j+1, I_{j+k}(x,y)=I_j(x,y), J_{j+k}(x,y)=J_j(x,y)$, for all $j\in \mathcal{A}, (x,y)\in\mathbf{R^2}$. \\

  Let $\mathbf{T^1}:=\mathbf{R}/\mathbf{Z}$. The main result in this paper are the following theorem.\\

  \noindent{\bf Theorem 2.1.} {\it Suppose that condition $({\rm \mathbf{H}})$ holds and that for each $n+1\leq i\leq 2n, p_i(t)\in C^{\gamma}(\mathbf{T^1})$ with $\gamma>1-\frac{1}{n}$, and for each $0\leq i\leq n, p_i(t)\in L(\mathbf{T^1})$. In addition, assume that \\

 \noindent{{\rm (i)}} \hs for $j=1,2,\cdots,k$, $I_j(x,y), J_j(x,y)\in C^5(\mathbf{R^2})$, and for some constant $E\gg 1$ and any non-negative integers $p, q${\rm (}$p+q\leq5${\rm )}, there is a constant $B=B(E)$ depending on $E$ such that
  $$ \left| \frac{\partial^{p+q} I_j(x,y)}{{\partial x^{p}}{\partial y^{q}}}\cdot (h_0(x,y))^{\frac{p}{2n+2}+\frac{q}{2}}\right|\leq B<+\infty,$$
   $$ \left|\frac{\partial^{p+q} J_j(x,y)}{{\partial x^{p}}{\partial y^{q}}}\cdot (h_0(x,y))^{\frac{p-n}{2n+2}+\frac{q}{2}}\right|\leq B<+\infty,$$
   whenever $x^2+y^2\geq E$, where $h_0(x,y)=\frac{1}{2(n+1)}x^{2n+2}+\frac{1}{2}y^2$.\\
 \noindent{{\rm (ii)}}\hs for $j=1,2,\cdots,k$, the impulsive maps $\Phi_{j}:(x,y)\mapsto (x,y)+(I_j(x,y),J_j(x,y))$
  are area-preserving homeomorphisms of $\mathbf{R^2}$.\\

    \noindent Then the time-{\rm 1} map $\widetilde{P}: (x, \dot{x})_{t=0}\mapsto (x, \dot{x})_{t=1}$ of Eq. {\rm (1.1)} possesses many {\rm (}positive Lebesgue measure{\rm )} invariant closed curves whose radiuses tend to infinity, and thus every solution $x(t)$ of {\rm (1.1)} is bounded for all time, i.e. it exists for all $t\in \mathbf{R}$ and
   $$\sup_{t\in \mathbf{R}}(|x(t)|+|\dot{x}(t)|)\leq \widetilde{C}<+\infty,  $$
where $\widetilde{C}=\widetilde{C}(x(0), \dot{x}(0))$ depends the initial data $(x(0), \dot{x}(0))$.}\\

        \noindent {\bf Remark 2.1.}  It is well-known that the solutions of autonomous unforced Duffing equations are a family of closed curves in the phase plane. Due to the existence of impulsive forces, the closed curves cannot be preserved usually. In [18], the choice of impulsive functions $I_j\equiv 0, J_j=-2\dot{x}(t_j^-)$ ensures that the trajectories of corresponding autonomous Duffing equation subject to impulsive forcing cannot escape the closed curve.  In [2], the trajectories of corresponding autonomous Duffing equation subject to impulsive forcing probably escape the closed curve, and Moser's twist theorem works for the equation studied. Compared with the those imposed on  the impulsive functions,  the conditions of the present Theorem 2.1  are  easier to be realized. Indeed,
 let $I_j(x,y)$ be independent of $(x,y)$, that is $I_j(x,y)=\alpha_j\in\mathbf{R}$, and let\\
            $$ J_j(x,y)=\sum_{m=0}^n\beta_{jm} x^m, \hs \beta_{jm}\in\mathbf{R}, $$\\
where $j=1,2,\cdots,k$. Then conditions (i) and (ii) of Theorem 2.1 are clearly satisfied. More examples can be given. For example, if $n\geq 5$, then the following impulsive functions:\\
$$  I_j(x,y)=\alpha_j,\hs J_j(x,y)=\beta_j\sin(x+\gamma_j)\hp {\rm or}\hp J_j(x,y)=\beta_je^{-x^2}, \beta_je^{-x^4},\cdots, $$\\
  where $\alpha_j,\beta_j,\gamma_j\in \mathbf{R}, j=1,2,\cdots,k$, satisfy conditions (i) and (ii) of Theorem 2.1.\\

        It is interesting to observe that the above impulsive functions $I_j$'s and $J_j$'s do not satisfy the condition (i) of Theorem 1.1 in [2]. \\

         \noindent {\bf Remark 2.2.} If $I_j(x,y)\equiv 0, J_j(x,y)\equiv 0, j=1, 2, \cdots, k$, then conditions (i) and (ii) of Theorem 2.1 naturally hold  and Eq. (1.1) is just the unforced Duffing equation (1.2). Therefore, Theorem 2.1 generalizes the result by Yuan [29] for (1.2) to the impulsive forcing analogous ones.\\

\noindent{\bf Remark 2.3.} In condition (ii) of Theorem 2.1, the assumption that $\Phi_{j}:(x,y)\mapsto (x,y)+(I_j(x,y),J_j(x,y))$ are area-preserving is equivalent to\\
  $$\frac{\partial I_j}{\partial x}+\frac{\partial J_j}{\partial y}+\frac{\partial I_j}{\partial x}\cdot \frac{\partial J_j}{\partial y}-\frac{\partial I_j}{\partial y}\cdot \frac{\partial J_j}{\partial x}=0\hs {\rm or}\hskip 0.1cm -2 $$\\
for any $j=1,2,\cdots,k$ and $(x,y)\in\mathbf{R^2}$.\\

     \noindent {\bf Remark 2.4.} The equation (1.2) is a Hamiltonian system with the Hamiltonian function\\
\[H=h_0(x,y)+R(x,y,t),\]\\
where\
\[R(x,y,t)=\sum_{i=0}^{2n}\frac{1}{i+1}p_i(t)x^{i+1}.\]\\
We point out that $H$ can be generalized to more general function of $(x,y,t)$, not necessarily to be a polynomial of $(x,y)$. Here it is essential to regard $R$ as a relatively small perturbation with respect to $h_0$:
\[\left|\partial_x^{m_1}\partial_{y}^{m_2}\, \left(h_0(x,y)^{-1}\, R(x,y,t)\right) \right|<\infty,\; x^2+y^2\to \infty.\]
See [8] for the detail. Without this relatively small condition imposed on the perturbation $R$, the stability might have been violated. In the present paper, the condition (i) of Theorem 2.1 implies that the impulses are also relatively small with respect to $h_0$. A relatively large impulse might violate the stability, too. So the condition (i) of Theorem 2.1 should be reasonable. \\

    \noindent{\bf 3. Properties of solution operators}\\

In this section, we will discuss and establish some properties of solution operators for impulsive differential equations.\\

     \noindent{\it {\rm 3.1}. Local properties of solutions for impulsive differential equations}\\

    We begin with some basic results on the impulsive differential equation\\
 $$
 \left\{
 \begin{array}{ll}
 \dot{u}=F(t,u), \hs t\neq t_j,\\
 \bigtriangleup u(t_j)=L_j(u(t_j)), \hs j\in\mathcal{A},\\
 \bigtriangleup u(t_j):= u(t_j^+)-u(t_j), \hs j\in\mathcal{A}
  \end{array}
 \right.
 \eqno (3.1)
 $$\\
 with the initial value condition\\
 $$  u(\tau^+)=u_0,  \eqno (3.2) $$\\
where $\tau\in \mathbf{R}, u_0\in \mathbf{R^m}, m\in \mathbf{N}$, and where $u(\tau^+)=u(\tau)$ if $\tau\notin\mathcal{T}$. Let $\mathbf{I}\subseteq\mathbf{R}$ be an open interval, and $\mathbf{G}\subseteq\mathbf{R^m}$, an open connected set. Suppose that\\

\noindent $({\rm \mathbf{H_1}})$ $F: \mathbf{I}\times \mathbf{G}\mapsto \mathbf{R^m}$ is continuous, locally Lipschitz in the second variable.

\noindent $({\rm \mathbf{H_2}})$ $\{t_j\}$ is a strictly increasing sequence of real numbers, and $|t_j|\to \infty$ as $|j|\to \infty$.

\noindent $({\rm \mathbf{H_3}})$ The impulsive functions $L_j: \mathbf{G}\to \mathbf{R^m}$ are continuous for all $j\in\mathcal{A}$. \\

 Let $\mathbf{J}\subseteq\mathbf{I}$ be an open interval. \\

  \noindent{\bf Definition 3.1.} A function $\varphi: \mathbf{J}\to \mathbf{R^m}$ belongs to the set $\mathcal{PC}(\mathbf{J,\mathbf{\mathcal{T}}})$ if\\

     \noindent (a) it is left continuous;\\
     \noindent (b) it is continuous, except, possibly, points of $\mathbf{\mathcal{T}}\cap\mathbf{J}$, where it has discontinuities of the first kind. \\

 The Definition 3.1 means that if $\varphi(t)\in \mathcal{PC}(\mathbf{J,\mathbf{\mathcal{T}}})$, then $\varphi(t_j^+)=\lim_{t\to t_j+0}\varphi(t)$ exists and $\varphi(t_j^-)=\varphi(t_j)$, where $\varphi(t_j^-)=\lim_{t\to t_j-0}\varphi(t), t_j\in\mathbf{J}\cap\mathbf{\mathcal{T}}$.\\

     \noindent{\bf Definition 3.2.} A function $\varphi: \mathbf{J}\to \mathbf{R^m}$ belongs to the set $\mathcal{PC}^1(\mathbf{J,\mathbf{\mathcal{T}}})$ if
   $\varphi(t),\dot{\varphi}(t)\in \mathcal{PC}(\mathbf{J,\mathbf{\mathcal{T}}})$, where the derivative at points of $\mathbf{J}\cap\mathbf{\mathcal{T}}$ is assumed to be the left derivative. \\

   Assume that $\psi: \mathbf{J}\to \mathbf{R^m}$ is a solution of (3.1).\\

    \noindent{\bf Lemma 3.1.} [1] {\it The solution $\psi$ of {\rm (3.1)} belongs to $\mathcal{PC}^1(\mathbf{J,\mathbf{\mathcal{T}}})$.}\\

\noindent{\bf Proposition 3.1.} [1] {\it A function $\psi(t)\in\mathcal{PC}^1(\mathbf{J,\mathbf{\mathcal{T}}}), \psi(\tau^+)=u_0$, is a solution of {\rm (3.1)} if and only if\\
$$ \psi(t)=
\left\{
\begin{array}{l}
u_0+\int_{\tau}^tF(s,\psi(s))ds+\sum_{\tau\leq t_j<t}L_j(\psi(t_j)),\hs t\geq \tau,\\
u_0+\int_{\tau}^tF(s,\psi(s))ds-\sum_{t\leq t_j<\tau}L_j(\psi(t_j)),\hs t<\tau.
\end{array}
\right.
$$}

  We have the following local properties of solutions for equation (3.1).\\

\noindent {\bf Lemma 3.2.} (Local properties) [1, 7, 15] {\it Assume that the conditions $({\rm \mathbf{H_1}})$-$({\rm \mathbf{H_3}})$ hold and that every jump equation\\
$$   u=v+L_j(v),   \eqno (3.3)   $$\\
$j\in\mathcal{A}, u\in\mathbf{G}$, has at most one solution with respect to $v\in \mathbf{G}$. Then\\

\noindent{{\rm (i)}} for any $\tau\in \mathbf{I}$ and $u_0\in \mathbf{G}$, there is a unique {\rm (}local{\rm )} solution $u=u(t;\tau,u_0)$ of Eq. {\rm (3.1)} satisfying the initial value condition {\rm (3.2)}, which exists on $(\tau-\alpha, \tau+\alpha)$ for some $\alpha>0$.\\
 \noindent{{\rm (ii)}} the map $Q_t: u_0\mapsto u(t;\tau,u_0)$ is continuous in $u_0$ for $t\in (\tau-\alpha, \tau+\alpha)\setminus\mathbf{\mathcal{T}}$.\\
\noindent{{\rm (iii)}} the solution has elastic property, that is, for any $b_0>0$, there is $r_{b_0}>0$ such that the inequality $|u_0|\geq r_{b_0}$ implies $|u(t;\tau,u_0)|\geq b_0$, for $t\in (\tau-\alpha, \tau+\alpha)\setminus\mathbf{\mathcal{T}}$.}\\

Subject to impulsive forcing differential equation (3.1), the equation\\
$$  \dot{z}=F(t,z)    \eqno (3.4)    $$\\
is called the corresponding continuous equation or the corresponding unforced differential equation, of (3.1).

It should be mentioned that the extension of a (local) solution of (3.1) over a maximal interval of existence is more complicated than that of (3.4). As an example, let us consider the following equation:\\
$$  \left\{
\begin{array}{l}
\dot{u}=1+u^2,\hs t\neq t_j:=\frac{j\pi}{4},\\
\bigtriangleup u(t)=-1,\hs t=\frac{j\pi}{4},\hs \forall j\in \mathcal{A},
\end{array}
\right.  \eqno (3.5)
$$\\
whose solution $u(t)$ satisfying $u(0)=0$ is continuable for all $t\in\mathbf{R}$. In fact, this solution has the form\\
$$  u(t)=\tan\left(t-\frac{j\pi}{4}\right),\hs t\in \left(\frac{j\pi}{4},\frac{(j+1)\pi}{4}\right],\hs j\in \mathcal{A}\cup\{0\}.  $$\\
However, the corresponding continuous equation \\
$$   \dot{z}=1+z^2,\hs t\in\mathbf{R}    \eqno (3.6)   $$\\
 has the solution $z(t)=\tan t$ satisfying $z(0)=0$, whose maximal interval of existence is $(-\pi/2,\pi/2)$. It is also interesting if we check the solution of (3.5) with $u(0)=u_0$, where $u_0\geq 1$. Indeed, since the corresponding solution of (3.6) with $z(0)=u_0$ is $z(t)=\tan(t+\beta)$, where $\beta=\arctan u_0$ such that $\pi/4\leq \beta<\pi/2$, whose maximal interval of existence is $(\beta-\pi/2, \pi/2-\beta)$ and so the corresponding integral curve cannot arrive at the first impulse moment $t_1=\pi/4$. As a result, the solution $u(t)$ of (3.5) with $u(0)=u_0\geq 1$ is not continuable to the moment $t_1$. Thus, the solution does not exist on $[\pi/4,\infty)$.

 For general cases, let $(\tau,u_0)\in \mathbf{I}\times\mathbf{G}$ be given. Denote by $u(t)=u(t;\tau,u_0)$ a solution of (3.1) satisfying $u(\tau^+)=u_0$. We
can extend the solution over a maximal interval of existence in the following way: \\

\noindent (1) $t$ is increasing, $t\geq \tau$;\\
  \noindent (2) $t$ is decreasing, $t\leq \tau$.\\

   The case (1), in its own turn, consists of two sub-cases:\\

   \noindent (1-1) $\tau\notin\mathbf{\mathcal{T}}$;

   \noindent (1-2) $\tau\in\mathbf{\mathcal{T}}, \tau=t_j$ for some $j\in\mathcal{A}$, say.\\

   In the sequel, we denote by $z(t;\kappa,\nu),\kappa\in\mathbf{I}, \nu\in\mathbf{G}$, a solution of (3.4) such that $z(\kappa;\kappa,\nu)=\nu$.

   Let us consider the sub-case (1-1). That is, the initial moment is not a discontinuous point. To be more concrete, suppose that $t_{j-1}<\tau<t_j$ for some fixed $j\in \mathcal{A}$. Denote by $[\tau,r),\tau<r$, the right maximal interval of the solution $z(t;\tau,u_0)$. If $r\leq t_j$, then $[\tau,r)$ is the maximal interval of $u(t)$, and $u(t)=z(t;\tau,u_0)$ for $t\in [\tau,r)$. Otherwise, $r>t_j$, and if $u(t_j)+L_j(u(t_j))\notin \mathbf{G}$, then $u(t)$ is not continuable beyond $t=t_j$, and the right maximal interval of $u(t)$ is $[\tau,t_j)$. If $u(t_j)+L_j(u(t_j))\in \mathbf{G}$, then $u(t)$ can be extended to the right as $u(t)=z(t;t_j,u(t_j^+))$, where $u(t_j^+)=u(t_j)+L_j(u(t_j))$. Denote by $[t_j,r)$ the right maximal interval of $z(t;t_j,u(t_j^+))$. If $r\leq t_{j+1}$, then $[\tau,r)$ is the right maximal interval of $u(t)$, and $u(t)=z(t;t_j,u(t_j^+))$ for $t\in [t_j,r)$. If $r>t_{j+1}$, and $u(t_{j+1})+L_{j+1}(u(t_{j+1}))\notin\mathbf{G}$, then the right maximal interval of $u(t)$ is $[\tau,t_{j+1})$. If $r>t_{j+1}$ and $u(t_{j+1})+L_{j+1}(u(t_{j+1}))\in\mathbf{G}$, then $u(t)$ can be continued as $z(t;t_{j+1},u(t_{j+1}^+))$, where $u(t_{j+1}^+)=u(t_{j+1})+L_{j+1}(u(t_{j+1}))$, and so on.

   Now, consider the sub-case (1-2), $\tau=t_j$ for some $j\in \mathcal{A}$. That is, the initial moment is a discontinuous point. Denote by $[\tau,r),\tau<r$, the right maximal interval of the solution $z(t;\tau,u_0)$. If $r\leq t_{j+1}$, then $[\tau^+,r)$ is the maximal interval of $u(t)$, and $u(t)=z(t;\tau,u_0)$ for $\tau<t<r$. If $r>t_{j+1}$, and if $u(t_{j+1})+L_{j+1}(u(t_{j+1}))\notin\mathbf{G}$, then $u(t)$ is not continuable beyond $t=t_{j+1}$, and the right maximal interval of $u(t)$ is $[\tau^+,t_{j+1}]$. If $u(t_{j+1})+L_{j+1}(u(t_{j+1}))\in\mathbf{G}$, then $u(t)$ can be extended to the right as $u(t)=z(t;t_{j+1},u(t_{j+1}^+))$, where $u(t_{j+1}^+)=u(t_{j+1})+L_{j+1}(u(t_{j+1}))$. The rest is as in the sub-case (1-1).

   Next consider the left extension of the solution. That is, the case (2). Fix $j\in\mathcal{A}$ such that $t_j<\tau\leq t_{j+1}$. Denote by $(l,\tau],l<\tau$, the left maximal interval of existence of the solution $z(t;\tau,u_0)$. If $l>t_j$, then $(l,\tau]$ is the left maximal interval of $u(t)$. If $l=t_j$ and $z(t_j^+;\tau,u_0)$ does not exist, then $(l,\tau]$ is the left maximal interval of $u(t)$. Otherwise, if $l=t_j$ and $z(t_j^+;\tau,u_0)$ exists, or if $l<t_j$, then $u(t_j^+)=z(t_j^+;\tau,u_0)$. Now, if the jump equation\\
      $$  u(t_j^+)=v+L_j(v)\eqno (3.7)   $$\\
     is solvable with respect to $v$ in $\mathbf{G}$, then take the solution of this equation as $u(t_j)$, and continue $z(t;t_j,u(t_j))$ of (3.4) at the left side of $t_j$. If (3.7) does not have a solution in $\mathbf{G}$, then $(t_j,\tau]$ is the left maximal interval of $u(t)$. Remark that (3.7) may have several (even infinitely many) solutions. Assume that one solution $v$ of (3.7) can be chosen, then proceeding in above way, one could finally construct a left maximal interval of $u(t)$.

     Unite the left and right maximal intervals to obtain the maximal interval of $u(t)$. It should be pointed out, as one can see from the above discussion, that the solution $u(t;\tau,u_0)$ of (3.1) may have more than one left and right maximal interval, and that the solution may have a right-closed maximal interval. \\

 \noindent{\it {\rm 3.2}. Global properties of solutions for impulsive differential equations}\\

 In order to explain why all the solutions of impulsive Duffing equation (1.1) exist for all $t\in \mathbf{R}$ under appropriate hypotheses on the impulsive functions, we first establish a general result for equation (3.1). We consider Eq. (3.4) under the following assumption:\\

 \noindent $({\rm \mathbf{H_4}})$ For any $(\tau, z_0)\in \mathbf{R}\times \mathbf{G}$, there exists a unique solution $z(t,\tau,z_0)$ of (3.4) with the maximal interval of existence $\mathcal{J}_a:=(\tau-a,\tau+a)$, such that $\inf\{a|(\tau,z_0)\in\mathbf{R}\times\mathbf{G}\}=\delta>0$, where $\delta$ is independent of $\tau$ and $z_0$.\\

  \noindent {\bf Lemma 3.3.} (Global properties) {\it Assume that the conditions $({\rm \mathbf{H_1}})$-$({\rm \mathbf{H_4}})$ hold with $\mathbf{I}=\mathbf{R}$. In addition, assume that\\

   \noindent{\rm (A)} $0<t_j-t_{j-1}<\delta$ for $\forall j\in \mathcal{A}$, where $t_0:=0$;\\
   \noindent{\rm (B)} $u\in\mathbf{G}$ implies that $u+L_j(u)\in\mathbf{G}$ for all $j\in \mathcal{A}$;\\
   \noindent{\rm (C)} for any $j\in\mathcal{A}$ and $u\in\mathbf{G}$, every jump equation {\rm (3.3)} has unique solution with respect to $v\in\mathbf{G}$.\\

   Then the following conclusions hold true:\\

   \noindent{\rm (a)} for any $\tau\in\mathbf{R}, u_0\in \mathbf{G}$, there is a unique {\rm (}global{\rm )} solution $u=u(t;\tau,u_0)$ of Eq. {\rm (3.1)} satisfying the initial value condition {\rm (3.2)}, and it exists for all $t\in \mathbf{R}$.\\
    \noindent{\rm (b)} the map $Q_t: u_0\mapsto u(t;\tau,u_0)$ is continuous in $u_0$ for $t\in \mathbf{R}\setminus\mathbf{\mathcal{T}}$.\\
    \noindent{\rm (c)} the solution has elastic property, that is, for any $b_0>0$, there is $r_{b_0}>0$ such that the inequality $|u_0|\geq r_{b_0}$ implies $|u(t;\tau,u_0)|\geq b_0$, for $t\in \mathbf{R}\setminus\mathbf{\mathcal{T}}$.}\\

 \noindent{\bf Proof.} We first prove (a) if:\\

 \noindent(1) $t$ is increasing, $t\geq \tau$;\\
 \noindent(2) $t$ is decreasing, $t\leq \tau$.\\

 For the case (1), let $\tau\in [t_{m-1},t_m)$ for some $m\in\mathcal{A}$, and let $z_0(t)$ be the solution of the initial value problem\\
  $$
 \left\{
 \begin{array}{l}
 \dot{z}=F(t,z),\hs t\geq \tau,\\
 z(\tau)=u_0,
 \end{array}
 \right.
 $$\\
 defined on $\mathcal{J}_{a_0}^R:=[\tau,\tau+a_0)(a_0>0)$. By conditions $({\rm \mathbf{H_4}})$ and (A), $t_m\in(\tau,\tau+a_0)$. Thus, $z_0(t_m)\in\mathbf{G}$. Using condition (B) yields $z_0(t_m)+L_m(z_0(t_m))\in\mathbf{G}$. Then we may let $z_1(t)$ be the solution of the second initial value problem\\
 $$
 \left\{
 \begin{array}{l}
  \dot{z}=F(t,z),\hs t\geq t_m,\\
  z(t_m)=z_0(t_m)+L_m(z_0(t_m)),
  \end{array}
  \right.
  $$\\
  defined on $\mathcal{J}_{a_1}^R:=[t_m,t_m+a_1)(a_1>0)$. It is clear that $t_{m+1}\in(t_m,t_m+a_1)$. Thus, $z_1(t_{m+1})\in\mathbf{G}$, and so $z_1(t_{m+1})+L_{m+1}(z_1(t_{m+1}))\in\mathbf{G}$. Then, we may again formulate an initial value problem of the form\\
   $$
 \left\{
 \begin{array}{l}
  \dot{z}=F(t,z),\hs t\geq t_{m+1},\\
  z(t_{m+1})=z_1(t_{m+1})+L_{m+1}(z_1(t_{m+1})),
  \end{array}
  \right.
  $$\\
whose solution $z_2(t)$ exists on some interval $\mathcal{J}_{a_2}^R:=[t_{m+1},t_{m+1}+a_2)(a_2>0)$. By repeating this procedure, by induction, we may obtain the solution $z_k(t)$ of the initial value problem\\
 $$
 \left\{
 \begin{array}{l}
  \dot{z}=F(t,z),\hs t\geq t_{m+n-1},\\
  z(t_{m+n-1})=z_{n-1}(t_{m+n-1})+L_{m+n-1}(z_{n-1}(t_{m+n-1})),
  \end{array}
  \right.
  $$\\
defined on $\mathcal{J}_{a_n}^R:=[t_{m+n-1},t_{m+n-1}+a_n)(a_n>0), n=1,2,\cdots$. Finally, we define $u(t)$ by\\
$$ u(t)=
 \left\{
 \begin{array}{lll}
  z_0(t),\hs \tau\leq t\leq t_m,\\
  z_1(t),\hs t_m<t\leq t_{m+1},\\
  \vdots\\
  z_n(t),\hs t_{m+n-1}<t\leq t_{m+n},\hs n=1,2,\cdots.
  \end{array}
  \right.
  $$\\
Then it is easy to verify that $u(t)$ is the solution of (3.1) defined on $[\tau,\infty)$, the right side of $\tau$, such that $u(\tau^+)=u_0$.

For the case (2), let $\tau\in [t_{m-1},t_m)$ for some $m\in\mathcal{A}$, and let $z_0(t)$ be the solution of the initial value problem\\
  $$
 \left\{
 \begin{array}{l}
 \dot{z}=F(t,z),\hs t\leq \tau,\\
 z(\tau)=u_0,
 \end{array}
 \right.
 $$\\
 defined on $\mathcal{J}_{a_0}^L:=(\tau-a_0,\tau](a_0>0)$. By conditions $({\rm \mathbf{H_4}})$ and (A), $t_{m-1}\in(\tau-a_0,\tau]$. Thus, $z_0(t_{m-1})\in\mathbf{G}$. Denote $u(t):=z_0(t)$ for $t_{m-1}\leq t\leq \tau$. From condition (C) we see that the jump equation\\
 $$  u(t_{m-1})=v+L_{m-1}(v)   $$\\
 has unique solution with respect to $v\in\mathbf{G}$, then take the solution of this equation as $u(t_{m-1}^-)$. Then we may let $z_1(t)$ be the solution of the second initial value problem\\
 $$
 \left\{
 \begin{array}{l}
  \dot{z}=F(t,z),\hs t\leq t_{m-1},\\
  z(t_{m-1})=u(t_{m-1}^-),
  \end{array}
  \right.
  $$\\
  defined on $\mathcal{J}_{a_1}^L:=(t_{m-1}-a_1,t_{m-1}](a_1>0)$. It is clear that $t_{m-2}\in(t_{m-1}-a_1,t_{m-1})$. Thus, $z_1(t_{m-2})\in\mathbf{G}$. Denote $u(t)=z_1(t)$ for $t_{m-2}\leq t\leq t_{m-1}$. From condition (C) we see that the jump equation\\
 $$  u(t_{m-2})=v+L_{m-2}(v)   $$\\
 has unique solution with respect to $v\in\mathbf{G}$, then take the solution of this equation as $u(t_{m-2}^-)$. Thus, we may again formulate an initial value problem of the form\\
   $$
 \left\{
 \begin{array}{l}
  \dot{z}=F(t,z),\hs t\leq t_{m-2},\\
  z(t_{m-2})=u(t_{m-2}^-),
  \end{array}
  \right.
  $$\\
whose solution $z_2(t)$ exists on some interval $\mathcal{J}_{a_2}^L:=(t_{m-2}-a_2, t_{m-2}](a_2>0)$. By repeating this procedure, by induction, we may obtain the solution $z_n(t)$ of the initial value problem\\
 $$
 \left\{
 \begin{array}{l}
  \dot{z}=F(t,z),\hs t\leq t_{m-n},\\
  z(t_{m-n})=u(t_{m-n}^-),
  \end{array}
  \right.
  $$\\
defined on $\mathcal{J}_{a_n}^L:=(t_{m-n}-a_n,t_{m-n}](a_n>0), n=1,2,\cdots$. Finally, we define $u(t)$ by\\
$$ u(t)=
 \left\{
 \begin{array}{lll}
  z_0(t),\hs t_{m-1}<t\leq \tau,\\
  z_1(t),\hs t_{m-2}<t\leq t_{m-1},\\
  \vdots\\
  z_n(t),\hs t_{m-n-1}<t\leq t_{m-n},\hs n=1,2,\cdots.
  \end{array}
  \right.
  $$\\
Then it is easy to verify that $u(t)$ is the solution of (3.1) defined on $(-\infty,\tau]$, the left side of $\tau$, such that $u(\tau^+)=u_0$.

Uniting the cases (1) and (2), we have that the solution $u=u(t;\tau,u_0)$ of (3.1) exists for $t\in \mathbf{R}$. Thus, the proof of conclusion (a) is complete.
The proof of conclusions (b) and (c) is similar to that of Lemma 3.2 as in [1, 7]. We omit the details.
 \qed\\

\noindent {\bf Remark 3.1.} It is interesting to observe that the example in (3.5) can be illustrated using Lemma 3.3.\\

The following corollary is a direct result of Lemma 3.3, which was used in the existing literature. For example, see [1, 7, 15].\\

\noindent {\bf Corollary 3.1.} {\it Under the conditions of Lemma {\rm 3.3}, except that $({\rm \mathbf{H_4}})$ and {\rm (A)} are replaced by the condition that all the solutions of {\rm (3.4)} exist for all $t\in\mathbf{R}$, the conclusions of Lemma {\rm 3.3} still hold true.}\\

\noindent {\bf Remark 3.2.} In equation (3.1), if $\dot{u}=F(t,u)$ is conservative and the impulsive maps $\widetilde{\Phi}_j: u\mapsto u+L_j(u)(j\in\mathcal{A})$ are  homeomorphisms of $\mathbf{R^m}$, then the map $Q_t$ in Lemmas 3.2-3.3 is a homeomorphism for $t\in \mathbf{R}\setminus\mathbf{\mathcal{T}}$.\\

\noindent{\it {\rm 3.3}. Global existence of solutions for impulsive Duffing equation}\\

We will deduce some global properties of impulsive Duffing equation (1.1) under condition $({\rm \mathbf{H}})$. To this end, by letting $y=\dot{x}$ and noting that $x(t_j^-)=x(t_j), y(t_j^-)=y(t_j)$, we can rewrite equation (1.1) as an equivalent system of the form\\
 $$
 \left\{
 \begin{array}{lll}
 \dot{x}=y,\hs t\neq t_j,\\
 \dot{y}=-x^{2n+1}-\sum_{i=0}^{2n}p_i(t)x^i,\hs t\neq t_j;\\
 \bigtriangleup x(t_j)=I_j(x(t_j), y(t_j)),\\
 \bigtriangleup y(t_j)=J_j(x(t_j),y(t_j)),\hs j=1,2,\cdots,k.
 \end{array}
 \right.
 \eqno (3.8)
 $$\\
\noindent {\bf Corollary 3.2.} {\it Suppose that condition $({\rm \mathbf{H}})$ holds and that for each $n+1\leq i\leq 2n, p_i(t)\in C^{\gamma}(\mathbf{T^1})$ with $\gamma>1-\frac{1}{n}$, and for each $0\leq i\leq n, p_i(t)\in L(\mathbf{T^1})$.  In addition, assume that for each $j=1,2,\cdots,k$ and any fixed $(u,v)\in \mathbf{R^2}$ the jump system\\
  $$
  \left\{
  \begin{array}{l}
  u=x+I_j(x,y),\\
  v=y+J_j(x,y)
  \end{array}
  \right.
  \eqno (3.9)
  $$\\
  has unique solution with respect to $(x,y)\in\mathbf{R^2}$. Then\\

   \noindent{\rm (a)} for any $\tau\in\mathbf{R}, (x_0,y_0)\in \mathbf{R^2}$, there is a unique solution $(x(t),y(t))=(x(t;\tau,x_0,\\y_0),y(t;\tau,x_0,y_0)$ of Eq. {\rm (3.8)} satisfying the initial condition $x(\tau^+)=x_0,y(\tau^+)=y_0$, and it exists for all $t\in \mathbf{R}$.\\
    \noindent{\rm (b)} the map $Q_t: (x_0,y_0)\mapsto (x(t;\tau,x_0,y_0),y(t;\tau,x_0,y_0)$ is continuous in $(x_0,y_0)$ for $t\in \mathbf{R}\setminus\mathbf{\mathcal{T}}$.\\
    \noindent{\rm (c)} the solution has elastic property, that is, for any $b_0>0$, there is $r_{b_0}>0$ such that the inequalities $|x_0|\geq r_{b_0},|y_0|\geq r_{b_0}$ implies that $|x(t;\tau,x_0,y_0)|\geq b_0$ and $|y(t;\tau,x_0,y_0)|\geq b_0$ for $t\in \mathbf{R}\setminus\mathbf{\mathcal{T}}$.}\\

  \noindent{\bf Proof.} By [29] we know that every solution $(x(t),y(t))$ of the unforced Duffing equation\\
  $$
 \left\{
 \begin{array}{lll}
 \dot{x}=y,\\
 \dot{y}=-x^{2n+1}-\sum_{i=0}^{2n}p_i(t)x^i
  \end{array}
 \right.
 \eqno (3.10)
 $$\\
 satisfying the initial value condition $(x(t_0),y(t_0))=(x_0,y_0)$, where $t_0\in \mathbf{R}, (x_0,y_0)\in\mathbf{R^2}$ is unique and exists for all $t\in\mathbf{R}$. Thus, by Lemma 3.3 or Corollary 3.1, we see that the conclusions of Corollary 3.2 hold true. \qed\\

\noindent {\bf Remark 3.3.} Under condition $({\rm \mathbf{H}})$, if for each $j=1,2,\cdots,k$, the impulsive maps $ \Phi_{j}:(x,y)\mapsto (x,y)+(I_j(x,y),J_j(x,y))$
  are homeomorphisms of $\mathbf{R^2}$, then for any fixed $(u,v)\in \mathbf{R^2}$ and each $j=1,2,\cdots,k$, jump system (3.9) has unique solution with respect to $(x,y)\in\mathbf{R^2}$.\\

  \noindent{\bf 4. Time-1 map of impulsive Duffing eqiuation}\\

 In this section, we will establish some properties of time-1 map for (3.8). To this end, we first describe the time-1 map of equation (3.1). Let us assume that conditions $({\rm \mathbf{H_1}})$-$({\rm \mathbf{H_3}})$ hold with $I=\mathbf{R},G=\mathbf{R^m}$ and that the following condition $({\rm \mathbf{H_5}})$ holds. \\

\noindent{$({\rm \mathbf{H_5}})$} $F$ is 1-periodic in the first variable, and there exists a positive integer $k$ such that  $0< t_1<t_2<\cdots<t_{k}<1, t_{j+k}=t_j+1$, and $L_{j+k}(u)=L_j(u)$ for $\forall j\in\mathcal{A}$ and $\forall u\in\mathbf{R^m}$.\\

 Let $u(t)=u(t;u_0)$ be the solution of (3.1) satisfying the initial condition $u(0)=u_0$. Let\\
 $$
 \left.
 \begin{array}{llll}
 P_0: u_0\mapsto u(t_1):=u_1,\\
 \widetilde{\Phi}_1: u_1\mapsto u_1+L_1(u_1)=u(t_1^+):=u_1^+,\\
 P_1: u_1^+\mapsto u(t_2):=u_2,\\
 \widetilde{\Phi}_2: u_2\mapsto u_2+L_2(u_2)=u(t_2^+):=u_2^+,\\
 \vdots\\
 P_{k-1}: u_{k-1}^+\mapsto u(t_k):=u_k,\\
 \widetilde{\Phi}_k: u_k\mapsto u_k+L_k(u_k)=u(t_k^+):=u_k^+,\\
 P_k: u_k^+\mapsto u(1).
 \end{array}
 \right.
 $$\\
 Then the time-1 map $P: u_0\mapsto u(1)$ of (3.1) can be expressed by\\
   $$ P=P_k\circ \widetilde{\Phi}_k\circ\cdots\circ P_1\circ\widetilde{\Phi}_1\circ P_0.   $$\\
     \noindent {\bf Remark 4.1.} Under condition $({\rm \mathbf{H_5}})$, if $\dot{u}=F(t,u)$ is conservative and the impulsive maps $\widetilde{\Phi}_j: u\mapsto u+L_j(u)(j=1,\cdots,k)$ are area-preserving homeomorphism of $\mathbf{R^m}$, then $P$ is an area-preserving homeomorphism of $\mathbf{R^m}$.\\

 In order to deduce the properties of the time-1 map of impulsive Duffing equation (3.8), we denote by $(x(t),y(t)=(x(t;x_0,y_0),y(t;x_0,y_0))$
 the solution of (3.8) satisfying the initial condition $(x(0),y(0)=(x_0,y_0)$. Let \\
 $$
 \left.
 \begin{array}{llll}
 \widetilde{P}_0: (x_0,y_0)\mapsto (x(t_1),y(t_1)):=(x_1,y_1),\\
 \Phi_1: (x_1,y_1)\mapsto (x_1+I_1(x_1,y_1),y_1+J_1(x_1,y_1))=(x(t_1^+),y(t_1^+)):=(x_1^+,y_1^+),\\
 \widetilde{P}_1: (x_1^+,y_1^+)\mapsto (x(t_2),y(t_2)):=(x_2,y_2),\\
 \Phi_2: (x_2,y_2)\mapsto (x_2+I_2(x_2,y_2),y_2+J_2(x_2,y_2))=(x(t_2^+),y(t_2^+)):=(x_2^+,y_2^+),\\
 \vdots\\
 \widetilde{P}_{k-1}: (x_{k-1}^+,y_{k-1}^+)\mapsto (x(t_k),y(t_k)):=(x_k,y_k),\\
 \Phi_k: (x_k,y_k)\mapsto (x_k+I_k(x_k,y_k),y_k+J_k(x_k,y_k))=(x(t_k^+),y(t_k^+):=(x_k^+,y_k^+),\\
 \widetilde{P}_k: (x_k^+,y_k^+)\mapsto (x(1),y(1)).
 \end{array}
 \right.
 $$\\
 Then the time-1 map $\widetilde{P}: (x_0,y_0)\mapsto (x(1),y(1))$ of (3.8) can be expressed by\\
  $$   \widetilde{P}=\widetilde{P}_k\circ \Phi_k\circ\cdots\circ \widetilde{P}_1\circ\Phi_1\circ \widetilde{P}_0.   $$\\
   \noindent{\bf Remark 4.2.} Under condition $({\rm \mathbf{H}})$, if the impulsive maps $ \Phi_{j}:(x,y)\mapsto (x,y)+(I_j(x,y),J_j(x,y))(j=1,2,\cdots,k)$ are area-preserving, then the time-1 map $\widetilde{P}$ of (3.8) is area-preserving, too.\\

From Corollary 3.2 and Remark 3.3, we have the following\\

\noindent {\bf Corollary 4.1.} {Suppose that the condition $({\rm \mathbf{H}})$ holds and that for each $n+1\leq i\leq 2n, p_i(t)\in C^{\gamma}(\mathbf{T^1})$ with $\gamma>1-\frac{1}{n}$, and for each $0\leq i\leq n, p_i(t)\in L(\mathbf{T^1})$. In addition, assume that the impulsive maps $ \Phi_{j}:(x,y)\mapsto (x,y)+(I_j(x,y),J_j(x,y))(j=1,2,\cdots,k)$ are homeomorphisms of $\mathbf{R^2}$. Then the time-1 map $\widetilde{P}$ of (3.8) is a homeomorphisms of $\mathbf{R^2}$. Moreover, for any $b_0>0$, there is $r_{b_0}>0$ such that the inequalities $|x_0|\geq r_{b_0},|y_0|\geq r_{b_0}$ implies that $|x(1;x_0,y_0)|\geq b_0$ and $|y(1;x_0,y_0)|\geq b_0$.}\\

 \noindent{\bf 5. Action-angle variables}\\

Let\
 $$x=AX.\eqno (5.1) $$\
 Then from equation (1.2), we get\\
  $$A\ddot{X}+A^{2n+1}X^{2n+1}+\sum_{i=0}^{2n} A^iX^{i}p_{i}(t)=0, $$
  here and in the sequel, $A$ is a constant large enough. That is,
  $$\ddot{X}+A^{2n}X^{2n+1}+\sum_{i=0}^{2n} A^{i-1}X^{i}p_{i}(t)=0. \eqno (5.2)$$\
Let\\
$$Y=A^{-n}\dot{X}=A^{-n-1}\dot{x}=A^{-n-1}y.\eqno (5.3)$$\\
Then, from (3.8), (5.1)-(5.3), we have\
\begin{eqnarray*}
\dot{Y}&=&A^{-n}\ddot{X}=A^{-n}(-A^{2n}X^{2n+1}-\sum_{i=0}^{2n} A^{i-1}X^{i}p_{i}(t))\\
&=&-A^{n}X^{2n+1}-\sum_{i=0}^{2n} A^{i-n-1}X^{i}p_{i}(t),
\end{eqnarray*}
$$\Delta X(t_{j}):=X(t_{j}^+)-X(t_{j})=A^{-1}\Delta x(t_{j})=A^{-1}I_j(AX(t_j),A^{n+1}Y(t_j)),$$
$$\Delta Y(t_{j}):=Y(t_{j}^+)-Y(t_{j})=A^{-n-1}\Delta y(t_j)=A^{-n-1}J_j(AX(t_j),A^{n+1}Y(t_j)).$$\
Thus,\\
$$\left\{\begin{array}{ll}
  \dot{X}=\frac{\partial H^*}{\partial Y}, \hs t\neq t_j,\\
  \dot{Y}=-\frac{\partial H^*}{\partial X}, \hs t\neq t_j,\\
  \Delta X(t_{j})=A^{-1}I_j(AX(t_j),A^{n+1}Y(t_j)):=\tilde{I}_j(X(t_j),Y(t_j)),\\
  \Delta Y(t_{j})=A^{-n-1}J_j(AX(t_j),A^{n+1}Y(t_j)):=\tilde{J}_j(X(t_j),Y(t_j)),
  \end{array}\right. \eqno (5.4) $$\\
where $j=1,2,\cdots,k$, and\\
$$ H^*(X,Y,t)=A^n\left(\frac{1}{2}Y^2+\frac{1}{2(n+1)}X^{2(n+1)}\right)+\sum_{i=0}^{2n}\frac{p_{i}(t)}{i+1} A^{i-n-1}X^{i+1}. \eqno (5.5)$$\\

\noindent{\bf Remark 5.1.} Let $\tilde{I}_j=\tilde{I}_j(X,Y),\tilde{J}_j=\tilde{J}_j(X,Y))$. Then, from (5.4), the condition (ii) of Theorem 2.1 and Remark 2.3, we have\\
 $$\frac{\partial \tilde{I}_j}{\partial X}+\frac{\partial \tilde{J}_j}{\partial Y}+\frac{\partial \tilde{I}_j}{\partial X}\cdot \frac{\partial \tilde{J}_j}{\partial Y}-\frac{\partial \tilde{I}_j}{\partial Y}\cdot \frac{\partial \tilde{J}_j}{\partial X}=0\hs {\rm or}\hskip 0.1cm -2.$$\\
 This implies that the absolute value $||\Delta||=1$ of Jacobian for impulsive maps $\widetilde{\Phi}_{j}^*:(X,Y)\mapsto (X,Y)+(\tilde{I}_j(X,Y),\tilde{J}_j(X,Y)), j=1,\cdots,k$. Thus each $\widetilde{\Phi}_{j}^*$ are area-preserving homeomorphisms of $\mathbf{R^2}$.\\

 Consider an auxiliary Hamiltonian system\\
  $$
  \left\{
  \begin{array}{l}
  \dot{X}=\frac{\partial H_0^*}{\partial Y},\\
   \dot{Y}=-\frac{\partial H_0^*}{\partial X},
   \end{array}
   \right.\eqno (5.6)
   $$\
  where\\
  $$H_0^*(X,Y)=\frac{1}{2}Y^2+\frac{1}{2(n+1)}X^{2(n+1)}.$$\\
  Let $(X_0(t), Y_0(t))$ be the solution of (5.6) with initial $(X_0(0), Y_0(0))=(1, 0)$. Then this solution
is clearly periodic. Let $T_0$ be its minimal positive period. By energy conservation, we have\\

\noindent $({\rm a_1})$  $(n+1)Y_0^{2}(t)+X_0^{2n+2}(t)\equiv 1$;\\
\noindent $({\rm a_2})$  $X_0(-t)=X_0(t)$,\hs $Y_0(-t)=-Y_0(t)$;\\
\noindent $({\rm a_3})$  $\dot{X}_0(t)=Y_0(t)$,\hs $\dot{Y}_0(t)=-X_0^{2n+1}(t)$;\\
\noindent $({\rm a_4})$  $X_0(t+T_0)=X_0(t)$,\hs $Y_0(t+T_0)=Y_0(t)$.\\

Let $\mathbf{T_{s}^{1}}=\left\{t\in \mathbf{C/Z}:|Im t|<s \right\}$ for any $s>0$.  We construct the following symplectic transformation\\
$$
  \va_0:
  \left\{\begin{array}{ll}
  X=c^{\alpha}\lambda^{\alpha}X_0(\ha T_0),\\
  Y=c^{\beta}\lambda^{\beta}Y_0(\ha T_0),
  \end{array}\right. \eqno (5.7) $$\\
where $\alpha=\frac{1}{n+2},\hs \beta=1-\alpha=\frac{n+1}{n+2},\hs c=\frac{1}{\alpha T_0}$ and where $(\lambda, \ha)\in \mathbf{R^+}\times \mathbf{T^1}$ is action-angle variables. By calculation, det$\frac{\partial(X,Y)}{\partial(\theta,\lambda)}=1$. Thus the transformation is indeed symplectic. Clearly $\Psi_0(\lambda,\ha)$ is analytic in $(\lambda, \ha)\in \mathbf{R^+}\times \mathbf{T^1_{s_0}}$ with some constant $s_0>1$.

By (5.7), we have\\
$$  \lambda=\frac{1}{c}[X^{2n+2}+(n+1)Y^2]^{\frac{n+2}{2n+2}}.  \eqno (5.8)  $$\\
We claim that there exists inverse function $\widetilde{X}_0^{-1}$ such that \\
$$   \theta=\widetilde{X}_0^{-1}(c^{-\alpha}\lambda^{-\alpha}X).     \eqno (5.9)   $$\\
Indeed, from (5.7) we have $X_0(\theta T_0)=c^{-\alpha}\lambda^{-\alpha}X$. We consider $\theta\in [0,1)$. For the case when $\theta\in [0,\frac{1}{2}]$, by ${\rm (a_3)}$ we get
$\frac{dX_0(\theta T_0)}{d\theta}=T_0Y_0(\theta T_0)<0$. Thus, we have $\theta=T_0^{-1}X_0^{-1}(c^{-\alpha}\lambda^{-\alpha}X)$. For the case when $\theta\in (\frac{1}{2},1)$, by using (5.7), ${\rm (a_2)}$ and ${\rm (a_4)}$ we have\\
$$   X=c^{\alpha}\lambda^{\alpha}X_0(\ha T_0)= c^{\alpha}\lambda^{\alpha}X_0(-\ha T_0)=c^{\alpha}\lambda^{\alpha}X_0((1-\ha)T_0). $$\\
Let $\xi=1-\theta$, then $\xi\in (0,\frac{1}{2})$. Thus, we have $\frac{dX_0(\xi T_0)}{d\xi}=T_0Y_0(\xi T_0)<0$ for $\xi\in (0,\frac{1}{2})$. Then we get\\
$$  \xi=T_0^{-1}X_0^{-1}(c^{-\alpha}\lambda^{-\alpha}X)\Longrightarrow \theta=1-T_0^{-1}X_0^{-1}(c^{-\alpha}\lambda^{-\alpha}X).  $$\

Now, from (5.4) we have that for $j=1,2,\cdots,k$\\
$$\left\{\begin{array}{ll}
  X(t_j^+)=X(t_j)+\tilde{I}_j(X(t_j),Y(t_j)),\\
  Y(t_j^+)=Y(t_j)+\tilde{J}_j(X(t_j),Y(t_j)).
  \end{array}\right.  \eqno (5.10) $$\

 Using (5.7) we have\\
$$
    \left\{\begin{array}{ll}
  X(t_j)=c^{\frac{1}{n+2}}\lambda^{\frac{1}{n+2}}(t_j)X_0(\ha(t_j)T_0),\\
  Y(t_j)=c^{\frac{n+1}{n+2}}\lambda^{\frac{n+1}{n+2}}(t_j)Y_0(\ha(t_j)T_0).
  \end{array}\right.  \eqno (5.11) $$\\
 Then using (5.7)-(5.11) we have that for $j=1,2,\cdots,k$\
 \begin{eqnarray*}
&&\bigtriangleup \lambda(t_j)=\lambda(t_j^+)-\lambda(t_j)\\
&=&\frac{1}{c}\{[X(t_j)+\tilde{I}_j(X(t_j),Y(t_j))]^{2n+2}+(n+1)[Y(t_j)+\tilde{J}_j(X(t_j),Y(t_j))]^2\}^{\frac{n+2}{2n+2}}\\
&&-\lambda(t_j)\\
&=&\frac{1}{c}\{[c^{\frac{1}{n+2}}\lambda^{\frac{1}{n+2}}(t_j)X_0(\ha(t_j)T_0)\\
&&+\tilde{I}_j(c^{\frac{1}{n+2}}\lambda^{\frac{1}{n+2}}(t_j)X_0(\theta(t_j)T_0),c^{\frac{n+1}{n+2}}
\lambda^{\frac{n+1}{n+2}}(t_j)Y_0(\theta(t_j)T_0))]^{2n+2}\\
&&+(n+1)[c^{\frac{n+1}{n+2}}
\lambda^{\frac{n+1}{n+2}}(t_j)Y_0(\ha(t_j)T_0)\\
&&+\tilde{J}_j(c^{\frac{1}{n+2}}\lambda^{\frac{1}{n+2}}(t_j)
X_0(\theta(t_j)T_0),c^{\frac{n+1}{n+2}}
\lambda^{\frac{n+1}{n+2}}(t_j)Y_0(\theta(t_j)T_0))]^2\}^{\frac{n+2}{2n+2}}-\lambda(t_j)\\
&:=&J^{*}_{j}(\lambda(t_j),\ha(t_j))
\end{eqnarray*}
and\\
\begin{eqnarray*}
&&\bigtriangleup\theta(t_j)=\theta(t_j^+)-\theta(t_j)\\
&=&\widetilde{X}_0^{-1}(c^{-\alpha}\lambda^{-\alpha}(t_j^+)X(t_j^+))-\theta(t_j)\\
&=&\widetilde{X}_0^{-1}(\{[c^{\frac{1}{n+2}}\lambda^{\frac{1}{n+2}}(t_j)X_0(\ha(t_j)T_0)\\
&&+\tilde{I}_j(c^{\frac{1}{n+2}}\lambda^{\frac{1}{n+2}}(t_j)X_0(\theta(t_j)T_0),c^{\frac{n+1}{n+2}}
\lambda^{\frac{n+1}{n+2}}(t_j)Y_0(\theta(t_j)T_0))]^{2n+2}\\
&&+(n+1)[c^{\frac{n+1}{n+2}}
\lambda^{\frac{n+1}{n+2}}(t_j)Y_0(\ha(t_j)T_0)\\
&&+\tilde{J}_j(c^{\frac{1}{n+2}}\lambda^{\frac{1}{n+2}}(t_j)
X_0(\theta(t_j)T_0),c^{\frac{n+1}{n+2}}
\lambda^{\frac{n+1}{n+2}}(t_j)Y_0(\theta(t_j)T_0))]^2\}^{-\frac{1}{2n+2}}\\
&&\cdot[c^{\frac{1}{n+2}}\lambda^{\frac{1}{n+2}}(t_j)X_0(\ha(t_j)T_0)\\
&&+\tilde{I}_j(c^{\frac{1}{n+2}}\lambda^{\frac{1}{n+2}}(t_j)X_0(\ha(t_j)T_0), c^{\frac{n+1}{n+2}}\lambda^{\frac{n+1}{n+2}}(t_j)Y_0(\ha(t_j)T_0))])\\
&&-\theta(t_j)\\
&:=&I^{*}_{j}(\lambda(t_j),\ha(t_j)).
\end{eqnarray*}

    Thus, under $\Psi_0$, equation (5.4) is changed into\\
  $$\left\{\begin{array}{ll}
  \dot{\ha}=\frac{\partial H}{\partial \lambda},\hs t\neq t_j,\\
  \dot{\lambda}=-\frac{\partial H}{\partial \ha},\hs t\neq t_j,\\
  \bigtriangleup \ha(t_j)=I^{*}_{j}(\lambda(t_j),\ha(t_j)),\\
  \bigtriangleup \lambda(t_j)=J^{*}_{j}(\lambda(t_j),\ha(t_j)),\hs j=1,2,\cdots,k,
  \end{array}\right.  \eqno (5.12)$$\\
where $H=H_0(\lambda)+R(\lambda,\ha,t)$ with\\
$$H_0(\lambda)=d\cdot A^n\cdot \lambda^{\frac{2(n+1)}{n+2}}, \hs d=\frac{c^{2\beta}}{2(n+1)} \eqno (5.13)$$
and\\
$$R(\lambda,\ha,t)=\sum_{i=0}^{2n}\frac{p_{i}(t)}{i+1} A^{i-n-1}(c^{\alpha}X_0(\ha T_0))^{i+1}\lambda^{\alpha(i+1)}.\eqno (5.14)$$\\
Clearly, $R(\lambda,\ha,t)=O(A^{n-1})$ for $A\rightarrow \infty$ and fixed $\lambda$ in some compact intervals.\\

\noindent {\bf 6. Approximation Lemma}\\

First, we cite an approximation lemma (see [22-23, 29] for the detail). We start by recalling
some definitions and setting some new notations. Assume that $X$ is a Banach space with the
norm $\|\cdot\|_X$. First recall that $C^{\mu}(\mathbf{R^n};X)$ for $0< \mu<1$ denotes the space of bounded H\"older
continuous functions $f:\mathbf{R^n}\rightarrow X$ with the form\\
$$\|f\|_{C^{\mu},X}=\sup_{0<|x-y|<1}\frac{\|f(x)-f(y)\|_X}{|x-y|^{\mu}}+\sup_{x\in \mathbf{R^n}}\|f(x)\|_X.$$\\
If $\mu=0$ then $\|f\|_{C^{\mu},X}$ denotes the sup-norm. For $\ell=m+\mu$ with $m\in \mathbf{N}$ and $0\leq\mu<1$, we denote
by $C^{\ell}(\mathbf{R^n},X)$ the space of functions $f:\mathbf{R^n}\rightarrow X$ with H\"older continuous partial derivatives,
i.e., $\partial^{\alpha}f\in C^{\mu}(\mathbf{R^n};X_{\alpha})$ for all multi-indices $\alpha=(\alpha_1,\cdots,\alpha_n)\in \mathbf{N^n}$ with the assumption that
$|\alpha|:=|\alpha_1|+\cdots+|\alpha_n|\leq m$, and $X_{\alpha}$ is the Banach space of bounded operators $T:\Pi^{|\alpha|}(\mathbf{R^n})\rightarrow X$ with the norm\\
 $$ \|T\|_{X_{\alpha}}=\sup\left\{\|T(u_1,u_2,\cdots,u_{|\alpha|})\|_X:\|u_i\|=1,1\leq i\leq|\alpha| \right\}. $$\\
  We define the norm $\|f\|_{C^{\ell}}=\sup_{|\alpha|\leq \ell}\|\partial^{\alpha}f\|_{C^{\mu},X_{\alpha}}$.\\

\noindent{\bf Lemma 6.1.} {\it {\rm( Jackson-Moser-Zehnder)} Let $f\in C^{\ell}(\mathbf{R^n};X)$ for some $\ell>0$ with finite $C^{\ell}$ norm over $\mathbf{R^n}$. Let $\phi$ be a radical-symmetric, $C^{\infty}$ function, having as support the closure of the
unit ball centered at the origin, where $\phi$ is completely flat and takes value {\rm 1}, and let $K=\hat{\phi}$ be its Fourier transform. For all $\sigma>0$ define\\
$$f_{\sigma}(x):=K_{\sigma}\ast f=\frac{1}{\sigma^n}\int_{\mathbf{T^n}}K(\frac{x-y}{\sigma})f(y)dy.$$\\
Then there exists a constant $C\geq 1$ depending only on $\ell$ and n such that the following holds: For any $\sigma>0$, the function $f_{\sigma}(x)$ is a real-analytic function from $\mathbf{C^n}$ to $X$ such that if $\Delta_{\sigma}^n$ denotes the n-dimensional complex strip of width $\sigma$,\\
$$\Delta_{\sigma}^n:=\left\{x\in \mathbf{C^n}||Imx_j|\leq\sigma,1\leq j\leq n \right\},$$\\
then for $\forall \alpha\in \mathbf{N^n}$ with $|\alpha|\leq \ell$ one has\\
$$\sup_{x\in \Delta_{\sigma}^n}\|\partial^{\alpha}f_{\sigma}(x)-\sum_{|\beta|\leq\ell-|\alpha|}\frac{\partial^{\beta+\alpha}
f(Rex)}{\beta!}(\sqrt{-1}Imx)^{\beta}\|_{X_{\alpha}}\leq C\|f\|_{C^{\ell}}\sigma^{\ell-|\alpha|}, $$\\
and for all $0\leq s\leq \sigma$,
$$\sup_{x\in \Delta_{s}^n}\|\partial^{\alpha}f_{\sigma}(x)-\partial^{\alpha}f_{s}(x)\|_{X_{\alpha}}\leq C\|f\|_{C^{\ell}}\sigma^{\ell-|\alpha|}. $$

The function $f_{\sigma}$ preserves periodicity {\rm (}i.e., if $f$ is T-periodic in any of its variable $x_j$, so is
$f_{\sigma}${\rm )}.}\\

By this lemma, for each $p_i\in C^{\gamma}(\mathbf{T^1})$, $i=n+1,n+2,\cdots,2n$, and any $\varepsilon>0$, there is a real analytic function (a complex value function $f(t)$ of complex variable $t$ in some domain in $\mathbf{C}$ is called real analytic if it is analytic in the domain and is real for real argument $t$) $p_{i,\varepsilon}(t)$ from $\mathbf{T_{\varepsilon}^1}$ to $\mathbf{C}$ such that
$$\sup_{t\in \mathbf{T^1}}|p_{i,\varepsilon}(t)-p_{i}(t)|\leq C\varepsilon^{\gamma}\|p_i\|_{C^{\gamma}} $$
and
$$\sup_{t\in \mathbf{T^1_{\varepsilon}}}|p_{i,\varepsilon}(t)|\leq C\|p_i\|_{C^{\gamma}}. $$

Write\\
$$R(\lambda,\ha,t)=R_{\varepsilon}(\lambda,\ha,t)+R^{\varepsilon}(\lambda,\ha,t), $$
where\\
$$R_{\varepsilon}(\lambda,\ha,t)=\sum_{i=n+1}^{2n}\frac{1}{i+1} A^{i-n-1}c^{\frac{i+1}{n+2}}X^{i+1}_0(\ha T_0)\lambda^{\frac{i+1}{n+2}}p_{i,\varepsilon}(t),$$
\begin{eqnarray*}R^{\varepsilon}(\lambda,\ha,t)&=&\sum_{i=0}^{n}\frac{1}{i+1} A^{i-n-1}c^{\frac{i+1}{n+2}}X^{i+1}_0(\ha T_0)\lambda^{\frac{i+1}{n+2}}p_{i}(t)\\
&&+\sum_{i=n+1}^{2n}\frac{1}{i+1} A^{i-n-1}c^{\frac{i+1}{n+2}}X^{i+1}_0(\ha T_0)\lambda^{\frac{i+1}{n+2}}(p_{i}(t)-p_{i,\varepsilon}(t)).
\end{eqnarray*}
Now let us restrict $\lambda$ to some compact intervals, [1,4], say. Let\\
$$  A^{-1}<\varepsilon_0.  \eqno (6.1)  $$\\
For a sufficiently small $\varepsilon_0>0$, letting\\
$$\varepsilon=\left(\frac{\varepsilon_0}{A^{n-1}}\right)^{\frac{1}{\gamma}}, $$\\
by Lemma 6.1, we have the following facts:\\

\noindent(I) $R^{\varepsilon}(\lambda,\ha,t)$ is real analytic in $(\lambda,\ha)\in [1,4]\times \mathbf{T^1_{s_0}}$ for fixed $t\in \mathbf{T^1}$ and $R^{\varepsilon}(\lambda,\ha,\cdot)\in L^1(\mathbf{T^1})$ for fixed $(\lambda,\ha)\in [1,4]\times \mathbf{T^1_{s_0}}$, and\\
$$\sup_{(\lambda,\ha,t)\in [1,4]\times \mathbf{T^1_{s_0}}\times \mathbf{T^1}}|R^{\varepsilon}(\lambda,\ha,t)|\leq C\varepsilon_0,   $$\\
where $C$ is a constant depending on only $\|p_i\|_{C^{\gamma}}$. Here and in the sequel, we denote by $C$ a universal constant which may be different in different place.\\

\noindent(II) $R_{\varepsilon}(\lambda,\ha,t)$ is real analytic in $(\lambda,\ha,T)\in [1,4]\times \mathbf{T^1_{s_0}}\times \mathbf{T^1_{\varepsilon}}$  and\\
$$\sup_{(\lambda,\ha,t)\in [1,4]\times \mathbf{T^1_{s_0}}\times \mathbf{T^1_{\varepsilon}}}|R_{\varepsilon}(\lambda,\ha,t)|\leq CA^{n-1}, $$\\
where $C$ depends on only $\|p_i\|_{C^{\gamma}}$. Therefore, we have\\
$$ H(\lambda,\ha,t)=H_0(\lambda)+R_{\varepsilon}(\lambda,\ha,t)+R^{\varepsilon}(\lambda,\ha,t).\eqno (6.2)$$\\
\noindent {\bf 7. Some symplectic transformations}\\

Following the idea of [29], we will look for a series of symplectic transformations $\Psi_1,\cdots,\Psi_N$ such that $H^{N}:=\Psi_N\circ\Psi_{N-1}\circ\cdots\circ
\Psi_1\circ H=H_0^N+O(\varepsilon_0)$, where $H_0^N(\mu)\approx A^n\mu^{\frac{2(n+1)}{n+2}}$ so that Moser's twist theorem works for $H^{N}$. To this end, we collect some conclusions verified in [29] as the following lemmas. For the detail, see Section 4 of [29].  \\

\noindent{\bf Lemma 7.1.} [29] {\it Let $ H(\lambda,\ha,t)$ be the same as {\rm (6.2)}. Then there is a symplectic diffeomorphism $\Psi_1$ depending periodically on $t$ of the form\\
$$
\Psi_1:
\left\{\begin{array}{ll}
  \lambda=\tilde{\mu}+u_1(\tilde{\mu},\tilde{\phi},t),\\
  \theta=\tilde{\phi}+v_1(\tilde{\mu},\tilde{\phi},t),
\end{array}\right.  $$\\
where $\sup_{D_1}|u_1|\leq CA^{-1}$, $\sup_{D_1}|v_1|\leq CA^{-1}$ and $D_1=[1+O(A^{-1}),4-O(A^{-1})]\times \mathbf{T^1_{s_0/2}}\times \mathbf{T^1_{\varepsilon/2}}$. Moreover the transformed Hamiltonian vectorfield $\Psi_1(X_{H})=X_{H^1}$ is of the form£º
$$H^1(\tilde{\mu},\tilde{\phi},t)=H^1_0(\tilde{\mu},t)+\widetilde{R}^1_{\varepsilon}(\tilde{\mu},\tilde{\phi},t)+\Psi_1\circ R^{\varepsilon},$$
where\\
$$\sup_{D_1}|\widetilde{R}^1_{\varepsilon}(\tilde{\mu},\tilde{\phi},t)|\leq C\varepsilon_0^{-\frac{1}{\gamma}}A^{n-1-\varpi},\hs \varpi:=n-\frac{n-1}{\gamma}.$$}\\

\noindent{\bf Lemma 7.2.} [29] {\it There are a series of symplectic transformations $\Psi_1,\cdots,\Psi_N${\rm (}$\Psi_{i}$ is similar to $\Psi_1$ of Lemma {\rm 7.1}, $i=2,\cdots,N${\rm )} with $n-\varpi N\leq -1,\hs N\in \mathbf{N}$, which transforms the Hamiltonian {\rm (6.2)} into
$$H^N(\mu,\phi,t)=\Psi_N\circ\cdots\circ\Psi_1\circ H=H_0^N(\mu,t)+R^N_{\varepsilon}(\mu,\phi,t)+\Psi_N\circ\cdots\circ\Psi_1\circ R^{\varepsilon},\eqno (7.1)$$
where \\
$$  H_0^N(\mu,t)=dA^n\mu^{\frac{2n+2}{n+2}}+O(A^{n-1}),  $$\\
and $A$ is large enough such that\\
$$A^{-1}\left(\frac{1}{\varepsilon_0}\right)^{\frac{N}{\gamma}}<\varepsilon_0.$$\\
Let\\
$$R^N(\mu,\phi,t)=R^N_{\varepsilon}(\mu,\phi,t)+\Psi_N\circ\cdots\circ\Psi_1\circ R^{\varepsilon},$$\\
then, we have\\
$$\sup_{(\mu,\phi,t)\in [2,3]\times \mathbf{T^1}\times \mathbf{T^1}}\int_0^1|\partial^p_\mu\partial^q_\phi R^N(\mu,\phi,t)|dt\leq C\varepsilon_0,\eqno (7.2)$$
\\
where $p,q$ are any non-negative integers with $0\leq p+q\leq 6$.}\\

   By Lemma 7.2, one easily see that, under the symplectic transformations $\Psi_1,\cdots,\Psi_N$, the corresponding unforced equation in (5.12) can be changed into\\
      $$\left\{\begin{array}{ll}
  \dot{\phi}=\frac{\partial H^N}{\partial \mu},\\
  \dot{\mu}=-\frac{\partial H^N}{\partial \phi},
   \end{array}
   \right.
   \eqno (7.3)
   $$\\
where $H^N(\mu,\phi,t)$ is as in (7.1).
Next we check the transformed impulsive forces in (5.12). First, from the symplectic transformation $\Psi_1$ in Lemma 7.1, by implicit function theorem we have\\
$$\left\{\begin{array}{ll}
  \tilde{\mu}=\lambda+u(\lambda,\theta,t),\\
  \tilde{\phi}=\theta+v(\lambda,\theta,t).
\end{array}\right.  $$\\
From this we have that, under the symplectic transformation $\Psi_1$, the jumps $\Delta\theta(t_j)$ and $\Delta\lambda(t_j)$ in (5.12) can be changed into\\
$$  \Delta\tilde{\phi}(t_j):=\tilde{\phi}(t_j^+)-\tilde{\phi}(t_j)=\tilde{I}_j^*(\tilde{\mu}(t_j),\tilde{\phi}(t_j)),   $$
 $$  \Delta\tilde{\mu}(t_j):=\tilde{\mu}(t_j^+)-\tilde{\mu}(t_j)=\tilde{J}_j^*(\tilde{\mu}(t_j),\tilde{\phi}(t_j)), $$\\
 where $j=1,2,\cdots,k$.

 In same way, under the symplectic transformation $\Psi_2$ which is similar to $\Psi_1$, the jumps $\Delta\tilde{\phi}(t_j)$ and $\Delta\tilde{\mu}(t_j)$ will be changed into new forms, say\\
 $$  \Delta\bar{{\phi}}(t_j):=\bar{{\phi}}(t_j^+)-\bar{{\phi}}(t_j)=\bar{{I}}_j^*(\bar{{\mu}}(t_j),\bar{{\phi}}(t_j)),$$
 $$   \Delta\bar{{\mu}}(t_j):=\bar{{\mu}}(t_j^+)-\bar{{\mu}}(t_j)=\bar{{J}}_j^*(\bar{{\mu}}(t_j),\bar{{\phi}}(t_j)), $$\\
where $j=1,2,\cdots,k$.

By repeating this procedure, we can see that, under the symplectic transformations $\Psi_1,\cdots,\Psi_N$ in Lemma 7.2, the jumps in (5.12) are finally changed into\\
  $$\left\{\begin{array}{ll}
    \bigtriangleup \phi(t_j):=\phi(t_j^+)-\phi(t_j)=I^{**}_{j}(\mu(t_j),\phi(t_j)),\\
  \bigtriangleup \mu(t_j):=\mu(t_j^+)-\mu(t_j)=J^{**}_{j}(\mu(t_j),\phi(t_j)),\\
\end{array}\right. \eqno (7.4)$$\\
where $j=1,2,\cdots,k$. Uniting (7.3) and (7.4), we see that, under $\Psi_1,\cdots,\Psi_N$, Eq. (5.12) can be transformed into
  $$\left\{\begin{array}{ll}
  \dot{\phi}=\frac{\partial H^N}{\partial \mu},\hs t\neq t_j,\\
  \dot{\mu}=-\frac{\partial H^N}{\partial \phi},\hs t\neq t_j,\\
  \bigtriangleup \phi(t_j)=I^{**}_{j}(\mu(t_j),\phi(t_j)),\\
  \bigtriangleup \mu(t_j)=J^{**}_{j}(\mu(t_j),\phi(t_j)),\hs j=1,2,\cdots,k\\
\end{array}\right. \eqno (7.5)$$

 It should be pointed out that, although we have not been able to formulated explicitly $I^{**}_{j}(\mu(t_j),\phi(t_j))$ and $J^{**}_{j}(\mu(t_j),\phi(t_j))$, we can implicitly express them. In order to make Moser's twist theorem works for the time-1 map of (7.5), in addition to the estimates of $H^N$ verified in Lemma 7.2, what we really need is establishing the estimates of the impulsive functions $I^{**}_{j}(\mu,\phi)$ and $J_j^{**}(\mu,\phi)$. We will treat it in detail in next section.\\

\noindent {\bf 8. Estimates of impulsive functions under symplectic transformations}\\

In this section, by using condition (i) of Theorem 2.1, we will establish some estimates for impulsive functions $I_j^{**}(\mu,\phi)$ and $J_j^{**}(\mu,\phi)$. To this end, we first give the estimates of $I_j^*(\lambda,\theta)$ and $J_j^*(\lambda,\theta)$. In this whole section and in the sequel, all the occurrences of $j$ mean $j=1,2,\cdots,k$.\\

\noindent{\bf Lemma 8.1.} {\it  For $(\lambda,\ha)\in [1,4]\times \mathbf{T^1}$, set \\
$$ \hat{I}_j(\lambda,\ha):=A^{-1}I_j(Ac^{\alpha}\lambda^{\alpha}X_0(\ha T_0),A^{n+1}c^{\beta}\lambda^{\beta}Y_0(\ha T_0)).$$
   $$\hat{J}_j(\lambda,\ha):=A^{-n-1}J_j(Ac^{\alpha}\lambda^{\alpha}X_0(\ha T_0),A^{n+1}c^{\beta}\lambda^{\beta}Y_0(\ha T_0)).$$\\
    Then for any non-negative integers $r, s${\rm (}$r+s\leq5${\rm )}, we have\\
    $$\left|\frac{\partial^{r+s} \hat{I}_j(\lambda,\ha)}{{\partial \lambda^{r}}{\partial \ha^{s}}}\right|,\hs\left|\frac{\partial^{r+s} \hat{J}_j(\lambda,\ha)}{{\partial \lambda^{r}}{\partial \ha^{s}}}\right|<CA^{-1},$$\\
if the condition {\rm (i)} of Theorem {\rm 2.1} holds.}\\

\noindent{\bf Proof.} For $(\lambda,\ha)\in [1,4]\times \mathbf{T^1}$, by (5.1), (5.3) and (5.7), and noting $(n+1)Y_0^{2}(t)+X_0^{2n+2}(t)\equiv 1$, we have\\
$$h_0(x,y)=\frac{1}{2(n+1)}(AX)^{2n+2}+\frac{1}{2}(A^{n+1}Y)^2=\frac{1}{2(n+1)}c^{\frac{2n+2}{n+2}}\lambda^{\frac{2n+2}{n+2}}A^{2n+2}.$$\\
Then from condition (i) of Theorem 2.1, for any non-negative integers $p$ and $q$($p+q\leq5$), we have\\
$$\left|\frac{\partial^{p+q} I_j(x,y)}{{\partial x^{p}}{\partial y^{q}}}\right|\leq CA^{-[p+(n+1)q]},\eqno (8.1)$$
   $$\left|\frac{\partial^{p+q} J_j(x,y)}{{\partial x^{p}}{\partial y^{q}}}\right|\leq CA^{-[p-n+(n+1)q]}.\eqno (8.2)$$\

    From (5.1), (5.3) and (5.7), we have\\
    $$  x=Ac^{\alpha}\lambda^{\alpha}X_0(\ha T_0),\hs y=A^{n+1}c^{\beta}\lambda^{\beta}Y_0(\ha T_0)).  $$ \

  Now applying the differential operator $D_\ha^sD_\lambda^r$ to $\hat{I}_j(\lambda,\ha)$ and $\hat{J}_j(\lambda,\ha)$, we consider four possible cases to prove the conclusion of Lemma 8.1.\\

   \noindent{\it Case {\rm 1:}}\hs $r=0,s=0$. From (8.1) and (8.2) we have\\
    $$ |\hat{I}_j(\lambda,\ha)|=|A^{-1}I_j(Ac^{\alpha}\lambda^{\alpha}X_0(\ha T_0),A^{n+1}c^{\beta}\lambda^{\beta}Y_0(\ha T_0))|=|A^{-1}I_j(x,y)|<CA^{-1},$$
   $$|\hat{J}_j(\lambda,\ha)|=|A^{-n-1}J_j(Ac^{\alpha}\lambda^{\alpha}X_0(\ha T_0),A^{n+1}c^{\beta}\lambda^{\beta}Y_0(\ha T_0))|=|A^{-n-1}J_j(x,y)|<CA^{-1}.$$\
  \noindent{\it Case {\rm 2:}}\hs $1\leq r\leq5,s=0$.  One can see that $\frac{\partial^r \hat{I}_j(\lambda,\ha)}{{\partial \lambda^r}}$ is a sum of terms
  $$ (D_y^qD_x^pA^{-1}I_j(x,y))(D_\lambda^{m_1}x)(D_\lambda^{m_2}x)\cdots(D_\lambda^{m_p}x)(D_\lambda^{n_1}y)(D_\lambda^{n_2}y)\cdots(D_\lambda^{n_q}y),$$
  where $1\leq p+q\leq r$ and $\sum_{i=1}^pm_i+\sum_{i=1}^qn_i=r$. From (8.1) we have\\
   $$\left|\frac{\partial^r \hat{I}_j(\lambda,\ha)}{{\partial \lambda^r}}\right|\leq CA^{-1-[p+(n+1)q]}\cdot A^{p}\cdot A^{(n+1)q}<CA^{-1}.$$\\
   Similarly, from (8.2) we have\\
   $$\left|\frac{\partial^r \hat{J}_j(\lambda,\ha)}{{\partial \lambda^r}}\right|<CA^{-1}.$$

  \noindent {\it Case {\rm 3:}} \hs $r=0,1\leq s\leq5$. One can see that $\frac{\partial^s \hat{I}_j(\lambda,\ha)}{{\partial \ha^s}}$ is a sum of terms\\
  $$ (D_y^qD_x^pA^{-1}I_j(x,y))(D_\ha^{m_1}x)(D_\ha^{m_2}x)\cdots(D_\ha^{m_p}x)(D_\ha^{n_1}y)(D_\ha^{n_2}y)\cdots(D_\ha^{n_q}y),$$\\
   where $1\leq p+q\leq s$ and $\sum_{i=1}^pm_i+\sum_{i=1}^qn_i=s$. From (8.1) we have\\
   $$ \left|\frac{\partial^s \hat{I}_j(\lambda,\ha)}{{\partial \ha^s}}\right|\leq  CA^{-1-[p+(n+1)q]}\cdot A^{p}\cdot A^{(n+1)q}<CA^{-1}.$$\\
   Similarly, from (8.2) we have
   $$ \left|\frac{\partial^s \hat{J}_j(\lambda,\ha)}{{\partial \ha^s}}\right|<CA^{-1}.$$

  \noindent {\it Case {\rm 4:}} \hs $r\geq1, s\geq1, r+s\leq5$. One can see that $\frac{\partial^{r+s} \hat{I}_j(\lambda,\ha)}{{\partial \lambda^{r}}{\partial \ha^{s}}}$ is a sum of terms\\
$$ D_\theta^s[(D_y^qD_x^pA^{-1}I_j(x,y))(D_\lambda^{m_1}x)(D_\lambda^{m_2}x)\cdots(D_\lambda^{m_p}x)(D_\lambda^{n_1}y)(D_\lambda^{n_2}y)\cdots(D_\lambda^{n_q}y)], $$\\
   where $1\leq p+q\leq r$ and $\sum_{i=1}^pm_i+\sum_{i=1}^qn_i=r$. It is easily seen that
\begin{eqnarray*}
&&D_\theta^s[(D_y^qD_x^pA^{-1}I_j(x,y))(D_\lambda^{m_1}x)\cdots(D_\lambda^{m_p}x)(D_\lambda^{n_1}y)\cdots(D_\lambda^{n_q}y)]\\
&&= \sum_{i=0}^sC_s^i(D_\theta^iD_y^qD_x^pA^{-1}I_j(x,y))\cdot D_\theta^{s-i}[(D_\lambda^{m_1}x)\cdots(D_\lambda^{m_p}x)(D_\lambda^{n_1}y)\cdots(D_\lambda^{n_q}y)],
\end{eqnarray*}
where\\

\noindent{\rm (1)} if $i=0$, $D_\theta^iD_y^qD_x^p[A^{-1}I_j(x,y)]=D_y^qD_x^p[A^{-1}I_j(x,y)]$;\\
\noindent{\rm (2)} if $i\geq 1$, $D_\theta^iD_y^qD_x^p[A^{-1}I_j(x,y)]$ is a sum of terms\\
$$ (D_y^{q+\tilde{q}}D_x^{p+\tilde{p}}A^{-1}I_j(x,y))(D_\ha^{\tilde{m}_1}x)(D_\ha^{\tilde{m}_2}x)\cdots(D_\ha^{\tilde{m}_{\tilde{p}}}x)(D_\ha^{\tilde{n}_1}y)
(D_\ha^{\tilde{n}_2}y)\cdots(D_\ha^{\tilde{n}_{\tilde{q}}}y),$$\\
   with $1\leq \tilde{p}+\tilde{q}\leq i$ and $\sum_{j=1}^{\tilde{p}}\tilde{m}_j+\sum_{j=1}^{\tilde{q}}\tilde{n}_j=i$. Thus, from (8.1) we have\\
     $$ \left|\frac{\partial^{r+s} \hat{I}_j(\lambda,\ha)}{{\partial \lambda^{r}}{\partial \ha^{s}}}\right|<CA^{-1}.$$\\
Similarly, we can get\\
$$ \left|\frac{\partial^{r+s} \hat{J}_j(\lambda,\ha)}{{\partial \lambda^{r}}{\partial \ha^{s}}}\right|<CA^{-1}.$$\\
   Then, Lemma 8.1 is proved. \qed\\

\noindent{\bf Lemma 8.2.} {\it Assume that $f_1,f_2,u,v\in C^5(\mathbf{R^+}\times\mathbf{T^1})$ and that for any non-negative integers $r, s${\rm (}$r+s\leq5${\rm )}, \\
 $$\left|\frac{\partial^{r+s} f_1(\lambda,\ha)}{{\partial \lambda^{r}}{\partial \ha^{s}}}\right|,\hs \left|\frac{\partial^{r+s} f_2(\lambda,\ha)}{{\partial \lambda^{r}}{\partial \ha^{s}}}\right|,\hs \left|\frac{\partial^{r+s} u(\lambda,\ha)}{{\partial \lambda^{r}}{\partial \ha^{s}}}\right|,\hs\left|\frac{\partial^{r+s} v(\lambda,\ha)}{{\partial \lambda^{r}}{\partial \ha^{s}}}\right|<CA^{-1} $$\\
for $(\lambda,\ha)\in [1+O(A^{-1}),4-O(A^{-1})]\times \mathbf{T^1}$, then for such $(\lambda,\ha)$, we have\\

\noindent {\rm (1)} $\left|\frac{\partial^{r+s} [f_1(\lambda,\ha)\cdot f_2(\lambda,\ha)]}{{\partial \lambda^{r}}{\partial \ha^{s}}}\right|<CA^{-1};$\\
\noindent {\rm (2)} $\left|\frac{\partial^{r+s} f_1(\lambda+u,\ha+v)}{{\partial \lambda^{r}}{\partial \ha^{s}}}\right|<CA^{-1}$.}\\

\noindent{\bf Proof.} Consider $(\lambda,\ha)\in [1+O(A^{-1}),4-O(A^{-1})]\times \mathbf{T^1}$.\\

\noindent (1)  Applying $D_\ha^sD_\lambda^r$ to $f_1\cdot f_2$, we have\\
$$D_\theta^{s}D_\lambda^{r}[f_1\cdot f_2]=\sum_{i_2=0}^{s}\sum_{i_1=0}^{r}C_{s}^{i_2}C_{r}^{i_1}(D_\theta^{i_2}D_\lambda^{i_1}f_1)
(D_\theta^{s-i_2}D_\lambda^{r-i_1}f_2).$$\\
Thus,\\
$$D_\ha^sD_\lambda^r[f_1(\lambda,\ha)\cdot f_2(\lambda,\ha)]<CA^{-1}.$$\\
\noindent (2) One can see that $D_\ha^sD_\lambda^rf_1(\lambda+u,\ha+v)$ is a sum of terms\\
$$ D_\theta^s[(D_\phi^qD_\mu^pf_1(\mu,\phi))(D_\lambda^{m_1}\mu)(D_\lambda^{m_2}\mu)\cdots(D_\lambda^{m_p}\mu)(D_\lambda^{n_1}\phi)
(D_\lambda^{n_2}\phi)\cdots(D_\lambda^{n_q}\phi)],$$\\
   where $\mu=\lambda+u$, $\phi=\theta+v$, $0\leq p+q\leq r$ and $\sum_{i=1}^pm_i+\sum_{i=1}^qn_i=r$. It is easily seen that\\
\begin{eqnarray*}
&&D_\theta^s[(D_\phi^qD_\mu^pf_1(\mu,\phi))(D_\lambda^{m_1}\mu)\cdots(D_\lambda^{m_p}\mu)(D_\lambda^{n_1}\phi)\cdots(D_\lambda^{n_q}\phi)]\\
&&= \sum_{i=0}^sC_s^i(D_\theta^iD_\phi^qD_\mu^pf_1(\mu,\phi))\cdot D_\theta^{s-i}[(D_\lambda^{m_1}\mu)\cdots(D_\lambda^{m_p}\mu)(D_\lambda^{n_1}\phi)\cdots(D_\lambda^{n_q}\phi)],
\end{eqnarray*}\\
where $D_\theta^iD_\phi^qD_\mu^pf_1(\mu,\phi)$ is a sum of terms\\
$$ (D_\phi^{q+\tilde{q}}D_\mu^{p+\tilde{p}}f_1(\mu,\phi))(D_\ha^{\tilde{m}_1}\mu)(D_\ha^{\tilde{m}_2}\mu)\cdots(D_\ha^{\tilde{m}_{\tilde{p}}}
\mu)(D_\ha^{\tilde{n}_1}\phi)
(D_\ha^{\tilde{n}_2}\phi)\cdots(D_\ha^{\tilde{n}_{\tilde{q}}}\phi),$$\\
   with $0\leq \tilde{p}+\tilde{q}\leq i$ and $\sum_{j=1}^{\tilde{p}}\tilde{m}_j+\sum_{j=1}^{\tilde{q}}\tilde{n}_j=i$. Thus, we have\\
     $$ \left|\frac{\partial^{r+s} f_1(\lambda+u,\ha+v)}{{\partial \lambda^{r}}{\partial \ha^{s}}}\right|<CA^{-1}.$$\\
 This completes the proof of Lemma 8.2. \qed\\

\noindent{\bf Lemma 8.3.} {\it  Suppose that condition {\rm (i)} of Theorem {\rm 2.1} holds. Set $\lambda(t_j)=\lambda, \ha(t_j)=\ha$, then for any non-negative integers $r, s${\rm (}$r+s\leq5${\rm )}, we have\\

 \noindent {\rm (1)} $\left|\frac{\partial^{r+s} I_j^{*}(\lambda,\ha)}{{\partial \lambda^{r}}{\partial \ha^{s}}}\right|<CA^{-1}$;\\
 \noindent {\rm (2)} $\left|\frac{\partial^{r+s} J_j^{*}(\lambda,\ha)}{{\partial \lambda^{r}}{\partial \ha^{s}}}\right|<CA^{-1}$\\

\noindent for $(\lambda,\ha)\in [1+O(A^{-1}),4-O(A^{-1})]\times \mathbf{T^1}$, where $I_j^{*}(\lambda(t_j),\ha(t_j))$ and $J_j^{*}(\lambda(t_j),,\ha(t_j))$ are as in {\rm (5.12)}.
}\\

\noindent{\bf Proof.}  In order to simplify the calculation, set\\
 $$  \hat{I}_j=\hat{I}_j(\lambda,\theta), \hs \hat{J}_j=\hat{J}_j(\lambda,\theta);   $$
 $$  {I}_j^*={I}_j^*(\lambda,\theta), \hs {J}_j^*={J}_j^*(\lambda,\theta),  $$\\
 where $\hat{I}_j(\lambda,\theta), \hat{J}_j(\lambda,\theta)$ are explicitly given in Lemma 8.1, and ${I}_j^*(\lambda,\theta), {J}_j^*(\lambda,\theta)$ are explicitly given in (5.12), respectively.

We first prove conclusion (2).  By (5.7), (5.10) and noting the impulsive functions in (5.4) we have
\begin{eqnarray*} J^{*}_{j}(\lambda,\ha)&=&\frac{1}{c}\{[c^\frac{1}{n+2}\lambda^\frac{1}{n+2}X_0(\ha T_0)+\hat{I}_j]^{2n+2}+(n+1)[c^\frac{n+1}{n+2}\lambda^\frac{n+1}{n+2}Y_0(\ha T_0)+\hat{J}_j]^2\}^{\frac{n+2}{2n+2}}\\
  &&-\lambda \\&=&\frac{1}{c}\{c^\frac{2n+2}{n+2}\lambda^\frac{2n+2}{n+2}X_0^{2n+2}(\ha T_0)+\sum_{i=1}^{2n+2}C_{2n+2}^{i}\hat{I}_{j}^{i}[c^\frac{1}{n+2}\lambda^\frac{1}{n+2}X_0(\ha T_0)]^{2n+2-i}\\
  &&+(n+1)[c^\frac{2n+2}{n+2}\lambda^\frac{2n+2}{n+2}Y_0^{2}(\ha T_0)+2\hat{J}_jc^\frac{n+1}{n+2}\lambda^\frac{n+1}{n+2}Y_0(\ha T_0)+\hat{J}_{j}^{2}]\}^{\frac{n+2}{2n+2}}
  -\lambda \\&=&\frac{1}{c}\{c^\frac{2n+2}{n+2}\lambda^\frac{2n+2}{n+2}+\sum_{i=1}^{2n+2}C_{2n+2}^{i}\hat{I}_{j}^{i}[c^\frac{1}{n+2}\lambda^\frac{1}{n+2} X_0(\ha T_0)]^{2n+2-i}\\&&+(n+1)[2\hat{J}_jc^\frac{n+1}{n+2}\lambda^\frac{n+1}{n+2}Y_0(\ha T_0)+\hat{J}_{j}^{2}]\}^{\frac{n+2}{2n+2}}
  -\lambda \\&=&\lambda\{1+\sum_{i=1}^{2n+2}C_{2n+2}^{i}\hat{I}_{j}^{i}c^\frac{-i}{n+2}\lambda^\frac{-i}{n+2}X_0^{2n+2-i}(\ha T_0)\\
  &&+(2n+2)\hat{J}_j c^\frac{-n-1}{n+2}\lambda^\frac{-n-1}{n+2}Y_0(\ha T_0)+(n+1)c^\frac{-2n-2}{n+2}\lambda^\frac{-2n-2}{n+2}\hat{J} _{j}^{2}\}^{\frac{n+2}{2n+2}}
  -\lambda.\end{eqnarray*}\\
      Set\\
  \begin{eqnarray*}
  f(\lambda,\theta)&=&\sum_{i=1}^{2n+2}C_{2n+2}^{i}\hat{I}_{j}^{i}c^\frac{-i}{n+2}\lambda^\frac{-i}{n+2}X_0^{2n+2-i}(\ha T_0)+(2n+2)\hat{J}_j c^\frac{-n-1}{n+2}\lambda^\frac{-n-1}{n+2}Y_0(\ha T_0)\\
  &&+(n+1)c^\frac{-2n-2}{n+2}\lambda^\frac{-2n-2}{n+2}\hat{J} _{j}^{2}.
  \end{eqnarray*}\\
  Then by Taylor's formula we have
  \begin{eqnarray*} J^{*}_{j}(\lambda,\ha)&=&\lambda[1+f(\lambda,\theta)]^{\frac{n+2}{2n+2}}
  -\lambda\\&=&\lambda\{1+\frac{n+2}{2n+2}[1+\xi f(\lambda,\theta)]^{\frac{-n}{2n+2}}\cdot f(\lambda,\theta)\}
  -\lambda\\&=&\frac{n+2}{2n+2}[1+\xi f(\lambda,\theta)]^{\frac{-n}{2n+2}}\cdot f(\lambda,\theta)\cdot\lambda,\end{eqnarray*}\\
  where $0<\xi<1$. For $(\lambda,\ha)\in [1+O(A^{-1}),4-O(A^{-1})]\times \mathbf{T^1}$, from Lemmas 8.1-8.2, it is easily seen that for any non-negative integers $r, s${\rm (}$r+s\leq5${\rm )}, we have\\
 $$\left|\frac{\partial^{r+s} f(\lambda,\theta)}{{\partial \lambda^{r}}{\partial \ha^{s}}}\right|<CA^{-1}.$$\\
 Thus\\
$$ \left|\frac{\partial^{r+s} J_j^{*}(\lambda,\ha)}{{\partial \lambda^{r}}{\partial \ha^{s}}}\right|<CA^{-1}.\eqno (8.3)$$\\

 Next we prove conclusion (1). Using (5.7), (5.10) and noting the impulsive functions in (5.4) we have\\
  $$\left\{\begin{array}{ll}
  \hat{I}_j=c^{\alpha}\cdot(\lambda+J^{*}_{j})^{\alpha}\cdot X_0((\ha+I^{*}_{j}) T_0)-c^{\alpha}\cdot\lambda^{\alpha}\cdot X_0(\ha T_0),\\
  \hat{J}_j=c^{\beta}\cdot(\lambda+J^{*}_{j})^{\beta}\cdot Y_0((\ha+I^{*}_{j}) T_0)-c^{\beta}\cdot\lambda^{\beta}\cdot Y_0(\ha T_0).
  \end{array}\right.  \eqno (8.4)$$\\
We consider two possible cases.\\

  \noindent{\it Case {\rm 1}}: Assume that\\
   $$|x_0(\ha T_0)|\leq\sqrt[2n+2]{\frac{1}{2}}, \hs \sqrt{\frac{1}{2n+2}}\leq |y_0(\ha T_0)|\leq \sqrt{\frac{1}{n+1}}.  $$\\
   Then, from (8.4) one has by Taylor's formula\\
  $$\hat{I}_j=\alpha c^{\alpha}\cdot (\lambda+\xi_1J^{*}_{j})^{\alpha-1}\cdot X_0((\ha+\xi_1I^{*}_{j})T_0)\cdot J^{*}_{j}+T_0 c^{\alpha}\cdot (\lambda+\xi_1J^{*}_{j})^{\alpha}\cdot Y_0((\ha+\xi_1I^{*}_{j})T_0)\cdot I^{*}_{j},\eqno (8.5)$$\\
  where $0<\xi_1<1$. It follows \\
  $$I^{*}_{j}=\frac{\hat{I}_j-\alpha c^{\alpha}\cdot (\lambda+\xi_1J^{*}_{j})^{\alpha-1}\cdot X_0((\ha+\xi_1I^{*}_{j})T_0)\cdot J^{*}_{j}}{T_0 c^{\alpha}\cdot (\lambda+\xi_1J^{*}_{j})^{\alpha}\cdot Y_0((\ha+\xi_1I^{*}_{j})T_0)}.$$\\
   Noting $\sqrt{\frac{1}{2n+2}}\leq |y_0(\ha T_0)|\leq \sqrt{\frac{1}{n+1}}$ and using (8.3), (8.4) and Lemma 8.1, it is easily seen that $|y_0((\ha+\xi_1I^{*}_{j})T_0)|>c_1>0$. Then from (8.3) and Lemma 8.1 we get\\
  $$|I^{*}_{j}|< CA^{-1}.\eqno (8.6)$$\\
  From (8.5) we have\\
  $$\frac{\hat{I}_j}{T_0 c^{\alpha}\cdot (\lambda+\xi_1J^{*}_{j})^{\alpha}}=X_0((\ha+\xi_1I^{*}_{j})T_0)\cdot\frac{\alpha\cdot J^{*}_{j}}{T_0\cdot(\lambda+\xi_1J^{*}_{j})}+Y_0((\ha+\xi_1I^{*}_{j})T_0)\cdot I^{*}_{j}.\eqno (8.7)$$\\
  Set\\
   $$ g_1(\lambda,\theta)=\frac{\hat{I}_j}{T_0 c^{\alpha}\cdot (\lambda+\xi_1J^{*}_{j})^{\alpha}},\hs g_2(\lambda,\theta)=\frac{\alpha\cdot J^{*}_{j}}{T_0\cdot(\lambda+\xi_1J^{*}_{j})}.  $$\\
   Then, from (8.3) and Lemma 8.1, for any non-negative integers $r, s${\rm (}$r+s\leq5${\rm )},  we obtain\\
$$\left|\frac{\partial^{r+s} g_1(\lambda,\ha)}{{\partial \lambda^{r}}{\partial \ha^{s}}}\right|,\hs\left|\frac{\partial^{r+s} g_2(\lambda,\ha)}{{\partial \lambda^{r}}{\partial \ha^{s}}}\right|<CA^{-1}.$$\\

We claim that for any non-negative integers $r, s${\rm (}$r+s\leq5${\rm )},\\
$$\left|\frac{\partial^{r+s} I^{*}_{j}(\lambda,\ha)}{{\partial \lambda^{r}}{\partial \ha^{s}}}\right|<CA^{-1}.\eqno (8.8)$$\\
 Since (8.6) implies that $|D_\theta^{s}D_\lambda^{r}I^{*}_{j}|<CA^{-1}$ holds true for $s=r=0$, thus inductively, assume that $|D_\theta^{s}D_\lambda^{r}I^{*}_{j}|<CA^{-1}$ holds for $r+s\leq n$ ($0\leq n\leq 4$), then for $r+s=n+1$, applying $D_\theta^{s}D_\lambda^{r}$ to (8.7) we have
\begin{eqnarray*}
D_\theta^{s}D_\lambda^{r}g_1&=&\sum_{i_2=0}^{s}\sum_{i_1=0}^{r}C_{s}^{i_2}C_{r}^{i_1}[D_\theta^{i_2}D_\lambda^{i_1}X_0((\ha+\xi_1I^{*}_{j})T_0)]
(D_\theta^{s-i_2}D_\lambda^{r-i_1}g_2)\\&&+\sum_{i_2=0}^{s}\sum_{i_1=0}^{r}C_{s}^{i_2}C_{r}^{i_1}[D_\theta^{i_2}D_\lambda^{i_1}Y_0(
(\ha+\xi_1I^{*}_{j})T_0)](D_\theta^{s-i_2}D_\lambda^{r-i_1}I^{*}_{j}),\hskip 10mm (8.9)\end{eqnarray*}\\
where $D_\theta^{i_2}D_\lambda^{i_1}X_0((\ha+\xi_1I^{*}_{j})T_0)$ is a sum of terms\\
$$(D_\eta^{p+q}X_0)(D_\ha^{\tilde{m}_1}\eta)(D_\ha^{\tilde{m}_2}\eta)\cdots(D_\ha^{\tilde{m}_q}\eta)\cdot D_\theta^{i_2-i}[(D_\lambda^{m_1}\eta)(D_\lambda^{m_2}\eta)\cdots(D_\lambda^{m_{p}}\eta)],$$\\
and $D_\theta^{i_2}D_\lambda^{i_1}Y_0((\ha+\xi_1I^{*}_{j})T_0)$ is a sum of terms\\
$$(D_\eta^{p+q}Y_0)(D_\ha^{\tilde{m}_1}\eta)(D_\ha^{\tilde{m}_2}\eta)\cdots(D_\ha^{\tilde{m}_q}\eta)\cdot D_\theta^{i_2-i}[(D_\lambda^{m_1}\eta)(D_\lambda^{m_2}\eta)\cdots(D_\lambda^{m_{p}}\eta)],$$\\
 with $\eta=(\ha+\xi_1I^{*}_{j})T_0$, $0\leq p\leq i_1$, $0\leq q\leq i$, $0\leq i\leq i_2$, $\sum_{k=1}^pm_k=i_1$ and $\sum_{k=1}^q\tilde{m}_k=i$. One can see that the right hand of (8.9) has three highest order terms:\\
  $$ Y_0((\ha+\xi_1I^{*}_{j})T_0)\cdot D_\theta^{s}D_\lambda^{r}I^{*}_{j}, $$
  $$  T_0\cdot\xi_1\cdot Y_0((\ha+\xi_1I^{*}_{j})T_0)\cdot g_2\cdot D_\theta^{s}D_\lambda^{r}I^{*}_{j},  $$
  $$ -T_0\cdot\xi_1\cdot X_0^{2n+1}((\ha+\xi_1I^{*}_{j})T_0)\cdot I^{*}_{j}\cdot D_\theta^{s}D_\lambda^{r}I^{*}_{j}.  $$\\
  These and (8.9) yield\\
\begin{eqnarray*}D_\theta^{s}D_\lambda^{r}I^{*}_{j}&=&[Y_0((\ha+\xi_1I^{*}_{j})T_0)+T_0\cdot\xi_1\cdot Y_0((\ha+\xi_1I^{*}_{j})T_0)\cdot g_2\\&&-T_0\cdot\xi_1\cdot X_0^{2n+1}((\ha+\xi_1I^{*}_{j})T_0)\cdot I^{*}_{j}]^{-1} \cdot\{D_\theta^{s}D_\lambda^{r}g_1\\&&-\sum_{i_2=0}^{s}\sum_{i_1=0}^{r}C_{s}^{i_2}C_{r}^{i_1}[D_\theta^{i_2}D_\lambda^{i_1}X_0(
(\ha+\xi_1I^{*}_{j})T_0)](D_\theta^{s-i_2}D_\lambda^{r-i_1}g_2)\\&&-\sum_{i_2=0}^{s}\sum_{i_1=0}^{r}C_{s}^{i_2}C_{r}^{i_1}
[D_\theta^{i_2}D_\lambda^{i_1}Y_0((\ha+\xi_1I^{*}_{j})T_0)](D_\theta^{s-i_2}D_\lambda^{r-i_1}I^{*}_{j})\\&&+[Y_0((\ha+\xi_1I^{*}_{j})T_0)+T_0\cdot\xi_1\cdot Y_0((\ha+\xi_1I^{*}_{j})T_0)\cdot g_2\\&&-T_0\cdot\xi_1\cdot X_0^{2n+1}((\ha+\xi_1I^{*}_{j})T_0)\cdot I^{*}_{j}]\cdot D_\theta^{s}D_\lambda^{r}I^{*}_{j}\},\end{eqnarray*}\\
this implies that\\
$$|D_\theta^{s}D_\lambda^{r}I^{*}_{j}|<CA^{-1}$$\\
holds for $r+s=n+1$. Thus claim (8.8) follows.\\

\noindent{\it Case {\rm 2}}: Assume that\\
 $$  \sqrt[2n+2]{\frac{1}{2}}\leq|x_0(\ha T_0)|\leq1,\hs |y_0(\ha T_0)|\leq\sqrt{\frac{1}{2n+2}}.  $$\\
 Then, from (8.4) one has by Taylor's formula
\begin{eqnarray*}
\hat{J}_j&=&\beta c^{\beta}(\lambda+\xi_2J^{*}_{j})^{\beta-1}Y_0((\ha+\xi_2I^{*}_{j})T_0)\cdot J^{*}_{j}\\
&&-T_0 c^{\beta}(\lambda+\xi_2J^{*}_{j})^{\beta}X_0^{2n+1}((\ha+\xi_2I^{*}_{j})T_0)\cdot I^{*}_{j},
\end{eqnarray*}
 where $0<\xi_2<1$. Similar to the statement as in the Case 1, we can prove (8.8).

         Then, from (8.3) and (8.8) we see that Lemma 8.3 is proved.  \qed\\

\noindent{\bf Lemma 8.4.} {\it  Assume that the condition {\rm (i)} of Theorem {\rm 2.1} holds. Let $\mu(t_j)=\mu$,  $\phi(t_j)=\phi$, then for any non-negative integers $r, s${\rm (}$r+s\leq5${\rm )}, \\
 $$\left|\frac{\partial^{r+s} I^{**}_{j}(\mu,\phi)}{{\partial \mu^{r}}{\partial \phi^{s}}}\right|,\hs\left|\frac{\partial^{r+s} J^{**}_{j}(\mu,\phi)}{{\partial \mu^{r}}{\partial \phi^{s}}}\right|<CA^{-1}  $$\\
 whenever $(\mu,\phi)\in [1+O(A^{-1}),4-O(A^{-1})]\times \mathbf{T^1}$.}\\

\noindent{\bf Proof.}  From Lemma 7.1, the symplectic diffeomorphism $\Psi_1$ is of the form\\
$$\left\{\begin{array}{ll}
  \lambda=\tilde{\mu}+u_1(\tilde{\mu},\tilde{\phi},t),\\
  \theta=\tilde{\phi}+v_1(\tilde{\mu},\tilde{\phi},t),
\end{array}\right.  \eqno (8.10)$$\\
 where $\sup_{D_1}|u_1|\leq CA^{-1}$, $\sup_{D_1}|v_1|\leq CA^{-1}$. By the implicit function theorem we have\\
$$\left\{\begin{array}{ll}
  \tilde{\mu}=\lambda+u(\lambda,\theta,t),\\
  \tilde{\phi}=\theta+v(\lambda,\theta,t),
\end{array}\right.  \eqno (8.11)$$\\
where $|u|\leq CA^{-1},|v|\leq CA^{-1}$ for $(\lambda,\ha)\in [1+O(A^{-1}),4-O(A^{-1})]\times \mathbf{T^1}$. Next we will prove that for $(\lambda,\ha)\in [1+O(A^{-1}),4-O(A^{-1})]\times \mathbf{T^1}$, $(\tilde{\mu},\tilde{\phi})\in [1+O(A^{-1}),4-O(A^{-1})]\times \mathbf{T^1}$ and any non-negative integers $r, s${\rm (}$r+s\leq6${\rm )}, the the following estimates hold true:\\
$$\left|\frac{\partial^{r+s} u}{{\partial \lambda^{r}}{\partial \theta^{s}}}\right|,\hs\left|\frac{\partial^{r+s} v}{{\partial \lambda^{r}}{\partial \theta^{s}}}\right|,\hs\left|\frac{\partial^{r+s} v_1}{{\partial \tilde{\mu}^{r}}{\partial \tilde{\phi}^{s}}}\right|,\hs\left|\frac{\partial^{r+s} u_1}{{\partial \tilde{\mu}^{r}}{\partial \tilde{\phi}^{s}}}\right|<CA^{-1}.  \eqno (8.12)$$\\
Indeed, we have from [29]\\
$$\left\{\begin{array}{ll}
  \lambda=\tilde{\mu}+\frac{\partial S_1}{\partial \theta}=\tilde{\mu}+\nu(\tilde{\mu},\theta,t),\\
  \tilde{\phi}=\theta+\frac{\partial S_1}{\partial \tilde{\mu}}=\theta+g(\tilde{\mu},\theta,t),
\end{array}\right.  \eqno (8.13)$$\\
where \\ $$S_1(\tilde{\mu},\theta,t)=\int_0^\theta\frac{[R_\varepsilon](\tilde{\mu},t)-R_\varepsilon(\tilde{\mu},\theta,t)}{\partial_{\tilde{\mu}}H_0
(\tilde{\mu})}d\theta, \hs [R_\varepsilon](\tilde{\mu},t)=\int_0^1R_\varepsilon(\tilde{\mu},\theta,t)d\theta.\eqno (8.14)$$\\
From (8.11) and (8.13) we have\\
$$u=-\nu(\lambda+u,\theta,t),\eqno (8.15)$$\\
and $|D_{\tilde{\mu}}\nu|\leq\frac{1}{2}$, so that $u$ is
uniquely determined by the contraction principle. Moreover, the implicit
function theorem implies that $u$ is $C^6$ with respect to $(\lambda,\ha)\in [1+O(A^{-1}),4-O(A^{-1})]\times \mathbf{T^1}$. We claim that for any non-negative integers $r, s${\rm (}$r+s\leq6${\rm )}, \\
$$\left|\frac{\partial^{r+s} u}{{\partial \lambda^{r}}{\partial \theta^{s}}}\right|<CA^{-1},\hs (\lambda,\ha)\in [1+O(A^{-1}),4-O(A^{-1})]\times \mathbf{T^1}.  \eqno (8.16)$$\\
Indeed, applying $(D_{\lambda})^n$ to the equation (8.15), the right hand side is a sum of terms\\
$$(D_{\tilde{\mu}}^p\nu)\cdot D_{\lambda}^{j_1}(\lambda+u)\cdot D_{\lambda}^{j_2}(\lambda+u)\cdots D_{\lambda}^{j_p}(\lambda+u), \eqno (8.17)$$\\
with $1\leq p\leq n$ and $\sum_{i=1}^pj_i=n$. The highest order term is the one
with $p=1$, namely $(D_{\tilde{\mu}}\nu)\cdot D_{\lambda}^nu$. Noting $|u|<CA^{-1}$ and assuming that for $j\leq n-1$ the estimates
$|D_{\lambda}^ju|<CA^{-1}$ hold true, then inductively, from (8.13)-(8.15) we can conclude that the same estimate holds true for $j=n$. In fact, from (8.13) and (8.14) we have\\
$$|D_{\tilde{\mu}}^p\nu|<CA^{-1},$$\\
so we obtain\\
$$|D_{\lambda}^nu|\leq \frac{1}{|(1-D_{\tilde{\mu}}\nu)|}\cdot c_2A^{-1}<CA^{-1}.$$\\
The estimates of $(D_{\theta})^{j}(D_{\lambda})^{i}u$ can be proven similarly. Thus, the claim (8.16) follows.

Similarly, for $(\lambda,\ha)\in [1+O(A^{-1}),4-O(A^{-1})]\times \mathbf{T^1}$ and $(\tilde{\mu},\tilde{\phi})\in [1+O(A^{-1}),4-O(A^{-1})]\times \mathbf{T^1}$, one concludes that \\
$$\left|\frac{\partial^{r+s} v}{{\partial \lambda^{r}}{\partial \theta^{s}}}\right|,\hs\left|\frac{\partial^{r+s} v_1}{{\partial \tilde{\mu}^{r}}{\partial \tilde{\phi}^{s}}}\right|,\hs\left|\frac{\partial^{r+s} u_1}{{\partial \tilde{\mu}^{r}}{\partial \tilde{\phi}^{s}}}\right|<CA^{-1},\hs 0\leq r+s\leq 6. $$

 From Section 7, we have known that, under symplectic transformation $\Psi_1$, equation (5.12) can be changed into a new form, say\\
 $$
 \left\{
 \begin{array}{lll}
 \dot{\tilde{\phi}}=\frac{\partial H^1}{\partial\tilde{\mu}},\hs t\neq t_j,\\
  \dot{\tilde{\mu}}=-\frac{\partial H^1}{\partial\tilde{\phi}},\hs t\neq t_j,\\
  \Delta\tilde{\phi}(t_j):=\tilde{\phi}(t_j^+)-\tilde{\phi}(t_j)=\tilde{I}_j^*(\tilde{\mu}(t_j),\tilde{\phi}(t_j)),\\
  \Delta\tilde{\mu}(t_j):=\tilde{\mu}(t_j^+)-\tilde{\mu}(t_j)=\tilde{J}_j^*(\tilde{\mu}(t_j),\tilde{\phi}(t_j)),
  \end{array}
  \right.
  \eqno (8.18)
  $$\\
 where Hamiltonian $H^1$ satisfies $\Psi_1(X_{H})=X_{H^1}$ and is as in Lemma 7.1. For $\tilde{I}_j^*$ and $\tilde{J}_j^*$, see Section 7. From (8.11) we have\\
 $$\left\{\begin{array}{ll}
  \tilde{I}_j^*(\tilde{\mu}(t_j),\tilde{\phi}(t_j))=\theta(t_{j}^+)+v(\lambda(t_{j}^+),\theta(t_{j}^+),t_j)-\theta(t_j)-v(\lambda(t_{j}),\theta(t_{j}),t_j),\\
 \tilde{J}_j^*(\tilde{\mu}(t_j),\tilde{\phi}(t_j))=\lambda(t_{j}^+)+u(\lambda(t_{j}^+),\theta(t_{j}^+),t_j)-\lambda(t_j)-u(\lambda(t_{j}),\theta(t_{j}),t_j).
\end{array}\right.  \eqno (8.19)$$\\
Let $\tilde{\mu}(t_j)=\tilde{\mu}_j$, $\tilde{\phi}(t_j)=\tilde{\phi}_j$, $\lambda(t_j)=\lambda_j$, $\theta(t_{j})=\theta_j$, $f(\lambda,\ha,t_j)=\bar{f}(\lambda,\ha)$, where $f(\lambda,\ha,t)$ is any function. Then from (8.18) and (8.19) we have by Taylor's formula\\
\begin{eqnarray*} \tilde{I}_j^*(\tilde{\mu}_j,\tilde{\phi}_j)&= &I^*_j(\lambda_j,\theta_j)+\bar{v}(\lambda_j+J^*_j(\lambda_j,\theta_j),\theta_j+I^*_j(\lambda_j,\theta_j))-\bar{v}(\lambda_j,\theta_j) \\&=&I^*_j(\lambda_j,\theta_j)+\bar{v}_\lambda(\lambda_j+\xi_3\cdot J^*_j(\lambda_j,\theta_j),\theta_j+\xi_3\cdot I^*_j(\lambda_j,\theta_j))\cdot J^*_j(\lambda_j,\theta_j)\\&&+\bar{v}_\theta(\lambda_j+\xi_3\cdot J^*_j(\lambda_j,\theta_j),\theta_j+\xi_3\cdot I^*_j(\lambda_j,\theta_j))\cdot I^*_j(\lambda_j,\theta_j)\\
&=&I^*_j(\tilde{\mu}_j+\bar{u}_1,\tilde{\phi}_j+\bar{v}_1)+\bar{v}_\lambda(\tilde{\mu}_j+\bar{u}_1+\xi_3\cdot J^*_j(\tilde{\mu}_j+\bar{u}_1,\tilde{\phi}_j+\bar{v}_1),\\&&\tilde{\phi}_j+\bar{v}_1+\xi_3\cdot I^*_j(\tilde{\mu}_j+\bar{u}_1,\tilde{\phi}_j+\bar{v}_1))\cdot J^*_j(\tilde{\mu}_j+\bar{u}_1,\tilde{\phi}_j+\bar{v}_1)\\&&+\bar{v}_\theta(\tilde{\mu}_j+\bar{u}_1+\xi_3\cdot J^*_j(\tilde{\mu}_j+\bar{u}_1,\tilde{\phi}_j+\bar{v}_1),\\&&\tilde{\phi}_j+\bar{v}_1+\xi_3\cdot I^*_j(\tilde{\mu}_j+\bar{u}_1,\tilde{\phi}_j+\bar{v}_1))\cdot I^*_j(\tilde{\mu}_j+\bar{u}_1,\tilde{\phi}_j+\bar{v}_1),\end{eqnarray*}
$  \hskip 137mm (8.20)   $\\

\noindent where $0<\xi_3<1$; and\\
\begin{eqnarray*} \tilde{J}_j^*(\tilde{\mu}_j,\tilde{\phi}_j)&= &J^*_j(\lambda_j,\theta_j)+\bar{u}(\lambda_j+J^*_j(\lambda_j,\theta_j),\theta_j+I^*_j(\lambda_j,\theta_j))-\bar{u}(\lambda_j,\theta_j) \\&=&J^*_j(\lambda_j,\theta_j)+\bar{u}_\lambda(\lambda_j+\xi_4\cdot J^*_j(\lambda_j,\theta_j),\theta_j+\xi_4\cdot I^*_j(\lambda_j,\theta_j))\cdot J^*_j(\lambda_j,\theta_j)\\&&+\bar{u}_\theta(\lambda_j+\xi_4\cdot J^*_j(\lambda_j,\theta_j),\theta_j+\xi_4\cdot I^*_j(\lambda_j,\theta_j))\cdot I^*_j(\lambda_j,\theta_j)\\
&=&J^*_j(\tilde{\mu}_j+\bar{u}_1,\tilde{\phi}_j+\bar{v}_1)+\bar{u}_\lambda(\tilde{\mu}_j+\bar{u}_1+\xi_4\cdot J^*_j(\tilde{\mu}_j+\bar{u}_1,\tilde{\phi}_j+\bar{v}_1),\\&&\tilde{\phi}_j+\bar{v}_1+\xi_4\cdot I^*_j(\tilde{\mu}_j+\bar{u}_1,\tilde{\phi}_j+\bar{v}_1))\cdot J^*_j(\tilde{\mu}_j+\bar{u}_1,\tilde{\phi}_j+\bar{v}_1)\\&&+\bar{u}_\theta(\tilde{\mu}_j+\bar{u}_1+\xi_4\cdot J^*_j(\tilde{\mu}_j+\bar{u}_1,\tilde{\phi}_j+\bar{v}_1),\\&&\tilde{\phi}_j+\bar{v}_1+\xi_3\cdot I^*_j(\tilde{\mu}_j+\bar{u}_1,\tilde{\phi}_j+\bar{v}_1))\cdot I^*_j(\tilde{\mu}_j+\bar{u}_1,\tilde{\phi}_j+\bar{v}_1),\end{eqnarray*}
$  \hskip 137mm (8.21)   $\\

\noindent where $0<\xi_4<1$.  By (8.12) and by letting $t=t_j$, for $(\lambda,\ha)\in [1+O(A^{-1}),4-O(A^{-1})]\times \mathbf{T^1}$, $(\tilde{\mu}_j,\tilde{\phi}_j)\in [1+O(A^{-1}),4-O(A^{-1})]\times \mathbf{T^1}$ and any non-negative integers $r, s${\rm (}$r+s\leq5${\rm )}, we have\\
$$\left|\frac{\partial^{r+s} \bar{v}_{\ha}}{{\partial \lambda^{r}}{\partial \theta^{s}}}\right|,\hs\left|\frac{\partial^{r+s} \bar{v}_{\lambda}}{{\partial \lambda^{r}}{\partial \theta^{s}}}\right|,\hs\left|\frac{\partial^{r+s} \bar{u}_{\ha}}{{\partial \lambda^{r}}{\partial \theta^{s}}}\right|,\hs\left|\frac{\partial^{r+s} \bar{u}_{\lambda}}{{\partial \lambda^{r}}{\partial \theta^{s}}}\right|,\hs\left|\frac{\partial^{r+s} \bar{v}_1}{{\partial \tilde{\mu}_j^{r}}{\partial \tilde{\phi}_j^{s}}}\right|,\hs\left|\frac{\partial^{r+s} \bar{u}_1}{{\partial \tilde{\mu}_j^{r}}{\partial \tilde{\phi}_j^{s}}}\right|<CA^{-1}.  $$

Now, by using these estimates and applying Lemma 8.2 and Lemma 8.3 to (8.20) and (8.21), we have that for any non-negative integers $r, s${\rm (}$r+s\leq5${\rm )}\\
$$ \left|\frac{\partial^{r+s} \tilde{I}_j^*}{{\partial \tilde{\mu}_j^{r}}{\partial \tilde{\phi}_j^{s}}}\right|,\hs\left|\frac{\partial^{r+s} \tilde{J}_j^*}{{\partial \tilde{\mu}_j^{r}}{\partial \tilde{\phi}_j^{s}}}\right|<CA^{-1}  \eqno (8.22) $$\\
for $(\tilde{\mu}_j,\tilde{\phi}_j)\in [1+O(A^{-1}),4-O(A^{-1})]\times \mathbf{T^1}$. Similarly, from Section 7, we have known that, under symplectic transformation $\Psi_2$, equation (8.18) can be changed into a new form, say\\
 $$
 \left\{
 \begin{array}{lll}
 \dot{\bar{{\phi}}}=\frac{\partial H^2}{\partial\bar{{\mu}}},\hs t\neq t_j,\\
  \dot{\bar{{\mu}}}=-\frac{\partial H^2}{\partial\bar{{\phi}}},\hs t\neq t_j,\\
  \Delta\bar{{\phi}}(t_j):=\bar{{\phi}}(t_j^+)-\bar{{\phi}}(t_j)=\bar{{I}}_j^*(\bar{{\mu}}(t_j),\bar{{\phi}}(t_j)),\\
  \Delta\bar{{\mu}}(t_j):=\bar{{\mu}}(t_j^+)-\bar{{\mu}}(t_j)=\bar{{J}}_j^*(\bar{{\mu}}(t_j),\bar{{\phi}}(t_j)),
  \end{array}
  \right.
  \eqno (8.23)
  $$\\
where Hamiltonian $H^2$ satisfies $\Psi_2(X_{H^1})=X_{H^2}$, and for $\bar{{I}}_j^*$ and $\bar{{J}}_j^*$, see Section 7. Then, similar to the statement of deducing (8.22), we have that for any non-negative integers $r, s${\rm (}$r+s\leq5${\rm )} the estimates\\
$$ \left|\frac{\partial^{r+s}\bar{{I}}_j^*}{{\partial\bar{{\mu}}_j^{r}}{\partial\bar{{\phi}}_j^{s}}}\right|,\hs\left|\frac{\partial^{r+s}\bar{{J}}_j^*}{{\partial \bar{{\mu}}_j^{r}}{\partial\bar{{\phi}}_j^{s}}}\right|<CA^{-1}  \eqno (8.24) $$\\
hold true for $(\bar{{\mu}}_j,\bar{{\phi}}_j)\in [1+O(A^{-1}),4-O(A^{-1})]\times \mathbf{T^1}$.

 Finally, by repeating this procedure, and noting the fact that $\Psi_{N}\circ\Psi_{N-1}\circ\cdot\cdot\cdot\circ\Psi_{1}$ transforms Eq. (5.12) into (7.5) via Lemma 7.2 , we have that for any non-negative integers $r, s${\rm (}$r+s\leq5${\rm )}, the following estimates hold true\\
 $$  \left|\frac{\partial^{r+s} I^{**}_{j}}{{\partial \mu^{r}}{\partial \phi^{s}}}\right|,\hs\left|\frac{\partial^{r+s} J^{**}_{j}}{{\partial \mu^{r}}{\partial \phi^{s}}}\right|<CA^{-1}    $$\\
for $(\mu,\phi)\in [1+O(A^{-1}),4-O(A^{-1})]\times \mathbf{T^1}$.
This completes the proof of Lemma 8.4. \qed\\

\noindent {\bf 9. Proof of Theorem 2.1} \\

In this section, we will prove Theorem 2.1. First, we prove the following\\

\noindent{\bf Lemma 9.1.} {\it If the condition {\rm (ii)} of Theorem {\rm 2.1} holds,
  then the time-{\rm 1} map $\Phi^1$ of Eq. {\rm (7.5)} is area-preserving. Moreover, the map $\Phi^1$ has the intersection property on $\Omega=\left\{(\mu,\phi)|2\leq\mu\leq 3,\phi\in \mathbf{T^1}\right\}$, i.e. if $\Gamma$ is an embedded circle in $\Omega$ homotopic to a circle $\mu=const$. in $\Omega$, then $\Phi^1(\Gamma)\bigcap \Gamma\neq\emptyset$.}\\

\noindent{\bf Proof.}  In view of the discussions in Section 4, the time-1 map $\Phi^1$ of Eq. (7.5) is\\
$$\Phi^1=P_{k}^*\circ\Phi^{*}_{k}\circ\cdot\cdot\cdot\circ P_1^*\circ \Phi^{*}_{1}\circ P_0^*,$$
where\\
$$
 \left.
 \begin{array}{llll}
 P_0^*: (\mu(0),\phi(0))\mapsto (\mu(t_1),\phi(t_1)):=(\mu_1,\phi_1),\\
 \Phi_1^*: (\mu_1,\phi_1)\mapsto (\mu_1+I_1^{**}(\mu_1,\phi_1),\phi_1+J_1^{**}(\mu_1,\phi_1))=(\mu(t_1^+),\phi(t_1^+)):=(\mu_1^+,\phi_1^+),\\
 P_1^*: (\mu_1^+,\phi_1^+)\mapsto (\mu(t_2),\phi(t_2)):=(\mu_2,\phi_2),\\
 \Phi_2^*: (\mu_2,\phi_2)\mapsto (\mu_2+I_2^{**}(\mu_2,\phi_2),\phi_2+J_2^{**}(\mu_2,\phi_2))=(\mu(t_2^+),\phi(t_2^+)):=(\mu_2^+,\phi_2^+),\\
 \vdots\\
 P_{k-1}^*: (\mu_{k-1}^+,\phi_{k-1}^+)\mapsto (\mu(t_k),\phi(t_k)):=(\mu_k,\phi_k),\\
 \Phi_k^*: (\mu_k,\phi_k)\mapsto (\mu_k+I_k^{**}(\mu_k,\phi_k),\phi_k+J_k^{**}(\mu_k,\phi_k))=(\mu(t_k^+),\phi(t_k^+):=(\mu_k^+,\phi_k^+),\\
 P_k^*: (\mu_k^+,\phi_k^+)\mapsto (\mu(1),\phi(1)).
 \end{array}
 \right.
 $$\\
 From Remark 5.1, we know that the impulsive maps $\widetilde{\Phi}_{j}^*:(X,Y)\mapsto (X,Y)+(\tilde{I}_j(X,Y),\tilde{J}_j(X,Y))$ are area-preserving. Since each $\Psi_i (i=0,1,\cdots,N$) is symplectic, it follows that, under the symplectic diffeomorphism $\Psi_{N}\circ\Psi_{N-1}\circ\cdot\cdot\cdot\circ\Psi_{1}\circ\va_0$, $\Phi^{*}_{j}$ are also area-preserving. Since each ${P}_{l}^*(l=0,1,2,\cdots,k)$ is the flow of the unforced Hamiltonian system\\
$$
\left\{
 \begin{array}{l}
 \dot{\phi}=\frac{\partial H^N}{\partial \mu},\\
  \dot{\mu}=-\frac{\partial H^N}{\partial \phi},
 \end{array}
  \right.
  $$\\
 where $H^N(\mu,\phi,t)=H_0^N(\mu,t)+R^N_{\varepsilon}(\mu,\phi,t)+\Psi_N\circ\cdots\circ\Psi_1\circ R^{\varepsilon}$, it follows that each $P_{l}^*$ is area-preserving. Hence, the map $\Phi^1$ is area-preserving. This implies that $\Phi^1$ has the intersection property. \qed\\

 Let $(\mu(t),\phi(t))=(\mu(t,\mu,\phi),\phi(t,\mu,\phi)$ be the solution of Eq. (7.5) with the initial value $(\mu(0),\phi(0))
 =(\mu,\phi)$. Set $\phi_1=\phi(1),\mu_1=\mu(1)$.\\

\noindent{\bf Lemma 9.2.} {\it Assume that conditions of Theorem {\rm 2.1} hold. Then the time-{\rm 1} map $\Phi^1$ of the flow $\Phi^t$ of Eq. {\rm (7.5)} is of the form\\
$$ \Phi^1:
\left\{
\begin{array}{l}
\phi_1=\phi+\alpha(\mu)+F(\mu,\phi),\\
\mu_1=\mu+G(\mu,\phi).
\end{array}
\right.
\eqno (9.1)
$$\\
Moreover, $\dot{\alpha}(\mu)>0$ and for any non-negative integers $r,s (r+s\leq 5)$\\
$$  \left|\frac{\partial^{r+s}F(\mu,\phi)}{\partial\mu^r\partial\phi^s}\right|,\hs \left|\frac{\partial^{r+s}G(\mu,\phi)}{\partial\mu^r\partial\phi^s}\right|
<C\varepsilon_0,   $$\\
where $(\mu,\phi)\in [2,3]\times \mathbf{T^1}$. }\\

 \noindent{\bf Proof.} From Lemma 7.2, impulsive Hamiltonian equation (7.5) is\\
 $$\left\{\begin{array}{ll}
  \dot{\mu}=-\frac{\partial H^N}{\partial \phi}=-\frac{\partial R^N(\mu,\phi,t)}{\partial \phi},\hs t\neq t_j\\
  \dot{\phi}=\frac{\partial H^N}{\partial \mu}=\frac{\partial H^N_0(\mu,t)}{\partial \mu}+\frac{\partial R^N(\mu,\phi,t)}{\partial \mu},\hs t\neq t_j\\
  \Delta \mu(t_j)=J^{**}_{j}(\mu(t_j),\phi(t_j)),\\
  \Delta \phi(t_j)=I^{**}_{j}(\mu(t_j),\phi(t_j)).\\
  \end{array}\right. \eqno (9.2) $$\\
By integral calculation and noting Proposition 3.1, we have that for $0\leq t\leq 1$ \\
$$
\left\{
\begin{array}{l}
\phi(t)=\phi+\alpha(\mu,t)+\widetilde{F}(\mu,\phi,t),\\
\mu(t)=\mu+\widetilde{G}(\mu,\phi,t),
\end{array}
\right. \eqno (9.3)
$$\\
where \\
$$ \alpha(\mu,t)=\int_0^t\frac{\partial H_0^N(\mu+\widetilde{G},s)}{\partial\mu}ds=\int_0^t\frac{\partial H_0^N(\mu,s)}{\partial\mu}ds+\int_0^t\int_0^1
\frac{\partial^2H_0^N}{\partial\mu^2}(\mu+\tau \widetilde{G},s)\widetilde{G}d\tau ds,  \eqno (9.4)  $$\\
$$ \widetilde{F}(\mu,\phi,t)=\int_0^t\frac{\partial R^N}{\partial\mu}(\mu+\widetilde{G},\phi+\alpha+\widetilde{F},s)ds+C(t),  \eqno (9.5)  $$\\
$$  \widetilde{G}(\mu,\phi,t)=-\int_0^t\frac{\partial R^N}{\partial\phi}(\mu+\widetilde{G},\phi+\alpha+\widetilde{F},s)ds+D(t), \eqno (9.6)   $$\\
with\\
$$ C(t)=
\left\{
\begin{array}{l}
0,\hs 0\leq t\leq t_1,\\
\sum_{j=1}^iI_j^{**}(\mu(t_j),\phi(t_j)),\hs t_i<t\leq t_{i+1},\hs i=1,\cdots,k-1,\\
\sum_{j=1}^kI_j^{**}(\mu(t_j),\phi(t_j)),\hs t_k<t\leq 1,
\end{array}
\right.
$$\\
and\\
$$ D(t)=
\left\{
\begin{array}{l}
0,\hs 0\leq t\leq t_1,\\
\sum_{j=1}^iJ_j^{**}(\mu(t_j),\phi(t_j)),\hs t_i<t\leq t_{i+1},\hs i=1,\cdots,k-1,\\
\sum_{j=1}^kJ_j^{**}(\mu(t_j),\phi(t_j)),\hs t_k<t\leq 1.
\end{array}
\right.
$$\\
In (9.3), we let $t=1$, and set $\alpha(\mu,1)=\alpha(\mu), \widetilde{F}(\mu,\phi,1)=F(\mu,\phi), \widetilde{G}(\mu,\phi,1)=G(\mu,\phi)$. Then (9.1) clearly holds true, and one easily verifies that for
$(\mu,\phi)\in [2,3]\times\mathbf{T^1}$ these equations have an unique solution in the space $|\widetilde{F}|,|\widetilde{G}|\leq C\varepsilon_0$ by using the contraction
principle. From Lemma 7.2 we have\\
 $$H^N_0=d\cdot A^n\cdot \mu^{\frac{2n+2}{n+2}}+O(A^{n-1}). $$\\
Thus, $\dot{\alpha}(\mu)\geq CA^n>0$. Moreover, $\widetilde{F}$ and $\widetilde{G}$ are $C^5$ with respective to $(\mu,\phi)$. Finally, by using Picard iteration and Gronwall's inequality, we can see that the estimates of $F(\mu,\phi)$ and $G(\mu,\phi)$ can inductively be verified from (9.5), (9.6) and in view of (7.2) and Lemma 8.4. \qed\\

 Now let us state Moser's twist theorem. Let $\mathcal D$ be an annulus:
 \[\mathcal {D}\;: \; a\le r\le b,\quad 0<a<b.\] For convenience, we introduce for a function $h\in C^l(\mathcal {D})$ the norm
 \[|h|_l=\sup_{\mathcal{D},\; m+n\le l}\left|\frac{\partial^{m+n}}{\partial r^m\partial \theta^n} \right|.\]

 \noindent{\bf Moser's Twist Theorem.}  {\it  Let $\alpha(r)\in C^l$ and $|\partial_r\, \alpha(r)|\ge \nu>0$ on the annulus $\mathcal{D}$ for some $l$ with $l\ge 5$, and $\varepsilon$ be a positive number.

 Then there exists a $\delta>0$ depending on $\varepsilon, l, \alpha(r)$, such that any area-preserving mapping\\
   \[M:\; \begin{array}{l} \theta_1=\theta+2\pi \alpha(r)+f(r,\theta) \\  r_1=r+g(r,\theta) \end{array}\]\\
 of $\mathcal D$ into $\mathbf{R^2}$ with $f,g\in C^l$ and
 \[|f|_l+|g|_l\le \nu\, \delta\]
 possesses an invariant curve of the form
 \[r=c+u(\xi),\quad \theta=\xi+v(\xi)\]
 in $\mathcal D$ where $u,\, v$ are continuously differentiable, of period $2\pi$ and satisfy
 \[|u|_1+|v|_1<\varepsilon,\]
 and $c$ is a constant in $(a,b)$.  Moreover, the induced mapping of this curve is given by\
 \[\xi\to \xi+\omega\]\
  where $\omega$ is incommensurable with $2\pi$, and satisfies infinitely many conditions
 \[\left|\frac{\omega}{2\pi}-\frac{p}{q}\right|\ge \gamma\, q^{-\tau}\]
 with some positive $\gamma,\tau$, for all integers $q>0$, $p$. In fact, each choice of $\omega$ in the range of $\alpha(r)$ and satisfying the above inequalities give rise to such an invariant curve.}\\

The Moser's twist theorem above can be found in pp. 50-54
of [15] (also see [21]).  It should be pointed out that the $\delta$ does not depend on $\nu$.  It should be also noted that the period $2\pi$ can be replaced by any period $T$. In addition, ``any area-preserving mapping" can be relaxed to ``any mapping which has intersection property ".  We are now in a position to prove Theorem 2.1. Let $\nu=C A^n$. From Lemma 9.1 and Lemma 9.2, by Moser's twist theorem, $\Phi^1$ has an invariant curve $\widetilde{\Gamma}$ in the annulus $[2,3]\times \mathbf{T^1}$. Since $A$ can be
arbitrarily large, it follows that the time-1 map of the original system has an invariant curve
$\widetilde{\Gamma}_A$ in the annulus $[2A+C,3A-C]\times \mathbf{T^1}$ with $C$ being a constant independent of $A$. Choosing
a sequence $A =A_m\rightarrow \infty$ as $m\rightarrow \infty$, we have that there are countable many invariant curves
$\widetilde{\Gamma}_{A_m}$, clustering at $\infty$. Therefore any solution of the original system is bounded. This completes
the proof of Theorem 2.1.\\

\noindent {\bf Remark 9.1.} Any solutions starting from the invariant curves $\widetilde{\Gamma}_{A_m}$ ($m=1,2,\cdots$) are quasi-periodic
with frequencies $(1,\omega_m)$ in time $t$, where $(1,\omega_m)$ satisfies Diophantine conditions and
$\omega_m>CA^n_m$. Actually, the frequencies can form a positive Lebesgue set in $\mathbf{R}$.\\

\vskip 0.5cm
\noindent{\bf Acknowledgments}\\

 This work is supported by the National Natural Science Foundation of China (No. 11571088, No. 11771093
 ) and the Zhejiang Provincial Natural Science Foundation of China (LY14A010024).\\

\noindent{\bf References} \\
  \newcounter{cankao}
\begin{list}
{[\arabic{cankao}]}{\usecounter{cankao}\itemsep=0cm} \small
\item M. Akhmet, Principles of Discontinuous Dynamical Systems, Springer-New York, 2010.
\item L. Chen, J. Shen, Invariant tori of impulsive Duffing-type equations via KAM technique,  ArXiv: 1927156 [math.DS] 21 Jun 2017.
\item
R. Dieckerhoff and E. Zehnder, Boundedness of solutions via the twist
theorem, Ann. Sc. Norm.Super. Pisa 14(1)(1987), 79-85.
 \item Y. Dong, Sublinear impulse effects and solvability of boundary value problems for differential equations with impulses, J. Math. Anal. Appl. 264 (2001), 32-48.
\item F. Jiang, J. Shen and Y. Zeng, Applications of the Poincar$\acute{e}$-Birkhoff theorem to impulsive Duffing equations at resonance, Nonlinear Anal. Real World Appl. 13 (2012), 1292-1305.
\item S. Laederich and M. Levi, Invariant curves and time-dependent potential, Ergodic. Theory Dynam. Systems 11 (1991), 365-378.
     \item V. Lakshmikantham, D. Bainov and P. Simeonov, Theory of Impulsive Differential Equations, World Scientific, Singapore, 1989.
  \item M. Levi, Quasiperiodic motions in superquadratic time periodic potentials. Comm. Math. Phys. 143(1) (1991), 43-83.
\item X. Li, B. Liu and Y. Sun, The large twist theorem and boundedness of solutions for polynomial potentials with $C^1$ time dependent coefficients, J. Differential Equations (to appear).
    \item J. Littlewood, Some problems in real and complex analysis, Heath, Lexington, Mass. 1968.
    \item B. Liu, Boundedness for solutions of nonlinear Hill's equations with periodic forcing terms via Moser's twist theorem, Journal of Differential Equations 79(1989), 304-315.
\item B. Liu, Boundedness for solutions of nonlinear periodic differential equations via Moser's twist theorem, Acta Math. Sinica (N.S.) 8(1992), 91-98.
  \item G. Morris, A case of boundedness of Littlewood's problem on oscillatory differential equations, Bull.
Austral. Math. Soc. 14 (1976), 71-93.

\item J. Moser, On invariant curves of aera-preserving mapping of annulus, Nachr. Akad. Wiss. Gottingen
Math. Phys. 2 (1962), 1-20.
\item J. Moser, Stable and Random Motion in Dynamic Systems, Ann. of Math. Studies, Princeton Uni. Press,
Princeton, NJ, 1973.
 \item J. Nieto, Basic theory for nonresonace impulsive periodic problems of first order, J. Math. Anal. Appl. 205 (1997), 423-433.
  \item J. Nieto and D. O'Regan, Variational approach to impulsive differential equations, Nonlinear Anal. Real World Appl. 10 (2009), 680-690.
  \item Y. Niu, X. Li, Boundedness of solutions in impulsive Duffing equations with polynomial potentials and $C^1$ time dependent coefficients. ArXiv: 1706.06460v1 [math.DS] 18 Jun 2017.
  \item J. Norris, Boundedness in periodically forced second order conservative systems, J. London Math. Soc. 45 (2) (1992), 97-112.
   \item D. Qian, L. Chen and X. Sun, Periodic solutions of superlinear impulsive differential equations: A geometric approach, J. Differential Equations 258 (2015), 3088-3106.
     \item H. Russman, Uber invariante Kurven differenzierbarer Abbildungen eines Kreisringes, Nachr. Akad. Wiss.
Gottingen, Math. Phys. (2) 1970, 67-105.
\item D. Salamon, The Kolmogorov-Arnold-Moser theorem, Mathematical Physics Electronic Journal 10(3)
(2004), 1-37.
\item D. Salamon, E. Zehnder, KAM theory in configuration space, Commun. Math. Helv. 64(1)(1989),
84-132.
\item J. Sun, H. Chen and J. Nieto, Infinitely many solutions for second-order Hamiltonian system with impulsive effects, Math. Comput. Modelling 54 (2011), 544-555
    \item Y. Wang, Unboundedness in a Duffing equation with polynomial potentials, J. Differential Equations 160(2)(2000), 467-479.
    \item X. Yuan, Invariant tori of Duffing-type equations, Advances in Math. (China), 24(1995), 375-376.
    \item X. Yuan, Invariant tori of Duffing-type equations, J. Differential Equations 142 (2) (1998), 231-262.
 \item X. Yuan, Lagrange stability for Duffing-type equations, J. Differential Equations 160 (1) (2000), 94-117.
\item X. Yuan, Boundedness of solutions for Duffing equation
with low regularity in time, Chinese Annals of Mathematics, Series B, 38(5)(2017), 1037-1046.

\end{list}
\end{document}